\documentclass[11pt,a4paper]{article}
\usepackage{amsfonts}
\usepackage{amsmath,amssymb}
\usepackage{mathrsfs}
\usepackage{hyperref}
\usepackage{graphicx}
\usepackage{float}

\newtheorem{thm}{Theorem}[section]
\newtheorem{cor}[thm]{Corollary}
\newtheorem{lem}[thm]{Lemma}
\newtheorem{prop}[thm]{Proposition}

\newtheorem{rem}[thm]{Remark}

\numberwithin{equation}{section}\allowdisplaybreaks

\def\leq{\leqslant}
\def\ge{\geqslant}

\def\leq{\leqslant}
\def\geq{\geqslant}

\begin{document}

\title{\large\bf  Navier-Stokes Equation in Super-Critical Spaces $E^s_{p,q}$}

\author{\footnotesize \bf Hans G. Feichtinger$^{\dag}$,  Karlheinz Gr\"ochenig$^{\dag}$, Kuijie Li$^{\ddag}$  and Baoxiang Wang$^{\$,}$\footnote{B. Wang is the corresponding author. The project is supported in part by NSFC, grant 11771024} \\
\footnotesize
\it $^\dag$Faculty of Mathematics,  University of Vienna , Oskar-Morgenstern-Platz 1, A-1090 Wien\\
\footnotesize \it $^\ddag$School of Mathematical Sciences, Fudan University, Shanghai, 200433, China \\
\footnotesize \it $^\$$LMAM, School of Mathematical Sciences, Peking University, Beijing 100871, PR of China \\
\footnotesize
\text{ Emails: hans.feichtinger@univie.ac.at {\rm (H.F.)},  karlheinz.groechenig@univie.ac.at {\rm (K.G.)}} \\ \footnotesize \text{likjie@fudan.edu.cn {\rm (K.L.)},  wbx@math.pku.edu.cn {\rm (B.W.)}}
}
\maketitle

\thispagestyle{empty}
\begin{abstract}
\noindent In this paper we develop a new way to study the global existence and uniqueness for the Navier-Stokes equation (NS) and consider the initial data in a class of modulation spaces $E^s_{p,q}$ with exponentially decaying weights $(s<0, \ 1<p,q<\infty)$ for which the norms are defined by
$$
\|f\|_{E^s_{p,q}} = \left(\sum_{k\in \mathbb{Z}^d} 2^{s|k|q}\|\mathscr{F}^{-1} \chi_{k+[0,1]^d}\mathscr{F} f\|^q_p \right)^{1/q}.
$$
The space $E^s_{p,q}$ is a rather rough function space and cannot be treated as a subspace of tempered distributions. For example, we have the embedding $H^{\sigma}\subset E^s_{2,1}$ for all $\sigma<0$ and $s<0$. It is known that $H^\sigma$ ($\sigma<d/2-1$) is a super-critical space of NS, it follows that $ E^s_{2,1}$ ($s<0$) is also super-critical for NS.
We show that NS has a unique global mild solution if the initial data belong to $E^s_{2,1}$ ($s<0$) and their Fourier transforms are supported in $ \mathbb{R}^d_I:= \{\xi\in \mathbb{R}^d: \ \xi_i \geq 0, \, i=1,...,d\}$.   Similar results hold for the initial data in $E^s_{r,1}$ with $2< r \leq d$. Our results imply that NS has a unique global solution if the initial value $u_0$ is in $L^2$ with ${\rm supp} \,  \widehat{u}_0 \, \subset  \mathbb{R}^d_I$.  \\


\end{abstract}

\section{Introduction}

In this paper, we mainly use the time-frequency method to study the Cauchy problem for the  incompressible
Navier-Stokes equation (NS):
\begin{align}
u_t -\Delta  u+u\cdot \nabla u +\nabla p=0, \ \  {\rm div}\, u=0, \ \  u(0,x)= u_0,  \label{NS1}
\end{align}
where   $u=(u_1,...,u_d)$ denotes the flow
velocity vector and $p$ describes the scalar pressure,
$u_0=(u_0^1,...,u_0^d)$ is an initial velocity with $\mathrm{div} \,
u_0=0$. $u_t=\partial u/\partial t,$ $ \Delta $ denotes
the  Laplacian.
It is easy to see that \eqref{NS1} can be rewritten as the following
equivalent form:
\begin{align}
  u_t -\Delta u+\mathbb{P}\ \textrm{div}(u\otimes u)=0, \ \
u(0, x)= u_0,
 \label{NS2}
\end{align}
where $\mathbb{P}=I-{\rm div}\, \Delta^{-1}\nabla$ is the matrix operator
projecting onto the divergence free vector fields,  $I$ is identity
matrix. It is known that NS has the following scaling invariance: If $u$ is a solution of NS, then the scaling function
\begin{align}
u_\lambda (t,x) := \lambda u(\lambda^2 t, \lambda x), \ \ p_\lambda (t,x) := \lambda^2 p(\lambda^2t, \lambda x), \ \lambda>0 \label{nsscaling}
\end{align}
 is also the solutions of NS with initial data $\lambda u_0(\lambda \ \cdot)$. Recall that a function space $X$ defined on $\mathbb{R}^d$ is said to be a  critical space if the norm of $u_\lambda$ in $X$  is invariant under the dilation \eqref{nsscaling}, namely $\|u_\lambda(t, \cdot)\|_X \sim \|u(t, \cdot)\|_X$ for all $\lambda>0$. It is known that NS has a class of critical spaces such as
$\dot H^{d/2-1}, \  L^d, \ \dot B^{-1+d/p}_{p,q}$ and $ BMO^{-1}$. The spaces $\dot H^s $ with $ s<d/2-1$,    $L^r$ with $r<d$ and $\dot B^{s}_{p,q}$ with $s<-1+ d/p$ are said to be super-critical spaces for NS.

On the other hand, there is a class of function spaces which may have variant scalings for different functions, say modulation spaces $M^s_{p,q}$; cf. Sugimoto and Tomita \cite{SuTo07} (see also \cite{HaWa14}).  It is well-known that the short time Fourier transform (STFT) plays a crucial role in the theory of time frequency analysis, the STFT equipped with the $L^p$ norm generates modulation spaces (cf. \cite{Fei83}). A frequency-discrete version of the STFT, so-called frequency uniform decomposition operator is a very useful tool in the study of nonlinear PDE (cf. \cite{ChWaWaWo18,WaHu07,Wa06}).  Let $\sigma$ be a smooth cut-off function adapted to the unit cube $[-1/2, 1/2]^d$ and $\sigma =0$ outside the cube $[-3/4, 3/4]^d$.
We writet $\sigma_k =\sigma(\cdot - k)$  and assume that
\begin{align} \label{sigmak}
\sum_{k\in \mathbb{Z}^d} \sigma_k (\xi) \equiv 1, \  \ \forall \; \
\xi \in
\mathbb{R}^d.
\end{align}
The frequency uniform decomposition operators are defined in the following way:
\begin{align} \label{FUD}
\Box_k := \mathscr{F}^{-1} \sigma_k \mathscr{F}, \quad k\in {\Bbb Z}^d.
\end{align}
Let $\mathscr{S}({\Bbb R}^d)$ be the Schwartz space and $\mathscr{S}'({\Bbb R}^d)$ be its dual space. Let $s\in \mathbb{R}, \  1\leq p,q\leq \infty$.  The modulation space $M^s_{p,q}({\Bbb R}^d)$ consists of all $f\in \mathscr{S}'({\Bbb R}^d)$ for which the following norm is finite:
\begin{align}
  \|f\|_{M^s_{p,q}} :=  \left\|\{\langle k\rangle^s  \Box_k f\}_{k\in \mathbb{Z}^d}  \right\|_{\ell^q(L^p)}. \label{defmod1}
\end{align}
Noticing that $\|f\|_{M^s_{p,q}} \sim \left\|\{(I-\Delta)^{s/2} \Box_k f\}_{k\in \mathbb{Z}^d}  \right\|_{\ell^q(L^p)}$, we can also define
\begin{align}
  \|f\|_{\dot M^s_{p,q}} :=  \left\|\{  (-\Delta)^{s/2} \Box_k f\}_{k\in \mathbb{Z}^d}  \right\|_{\ell^q(L^p)}, \label{defmod1a}
\end{align}
which is said to be the hybrid Riesz potential-modulation spaces.
In this paper, we are mainly interested in the modulation spaces which have the exponential regularity, i.e., the spatial regularity $\langle k\rangle^s$ is replaced by $2^{s|k|}$ in the definition of $M^s_{p,q}$ in \eqref{defmod1} and we denote
\begin{align}
  \|f\|_{E^s_{p,q}} :=  \left\|\{ 2^{s|k|}  \Box_k f\}_{k\in \mathbb{Z}^d}  \right\|_{\ell^q(L^p)}. \label{defmod2}
\end{align}
In the case $s>0$, $E^s_{p,q}$ can be treated as a subspace of $\mathscr{S}'({\Bbb R}^d)$ so that the norm $\|\cdot\|_{E^s_{p,q}}$ is complete (cf. \cite{Wa06}).  However, if $s<0$, $E^s_{p,q}$ cannot be regarded as a subspace of $\mathscr{S}'({\Bbb R}^d)$, since the function like $\widehat{f} =2^{-s|\xi|/2}$ has the finite norm in $E^s_{p,q}$, but is not in $\mathscr{S}'({\Bbb R}^d)$.  Roughly speaking, in the case $s<0$, we can treat $E^s_{p,q}$ as the subspaces of the Gelfand-Shilov space, which contains the distributions with exponential growth so that $\Box_k f$ in \eqref{defmod2} is meaningful and $\|\cdot\|_{E^s_{p,q}}$ becomes a complete norm. In this paper we will give a reasonable definition of $E^s_{p,q}$ in the case $s<0$, which is eventually an equivalent characterization of modulation space with negative exponential weight, see \eqref{espq} in Section \ref{Emodspace}.

 Our goal is the study of NS with a class of initial data in $E^s_{p,1}$, $s<0$. It is worth to mention that $D^\alpha\delta \in E^s_{p,q}$ for any $\alpha\in \mathbb{Z}^d_+$ and $s<0$, where $\delta$ is the Dirac measure.
 It is known that \eqref{NS1} is equivalent to the following integral equation:
\begin{align} \label{NS3}
u(t) =H (t) u_0 +  \mathscr{A} \mathbb{P}\ {\rm div}\,(u\otimes u),
\end{align}
where
\begin{align} \label{NSnote1}
H(t):=e^{ t \Delta }=\mathscr{F}^{-1}e^{-t|\xi|^{2 }}\mathscr{F}, \quad   (\mathscr{A} f)(t):= \int_0^t H (t-\tau)f(\tau)d\tau.
\end{align}
The solutions of \eqref{NS3} are usually said to be mild solutions.  NS has been extensively studied  in \cite{BaBiTa12, Can97, CaPl96, Che99, DoLi09, EsSeSv03, FoTe89, GePaSt07, GrKu98, Iw10, Ka84, KoTa01, LeRi00, Pl96}. Cannone \cite{Can97} and Planchon \cite{Pl96} considered global solutions in 3D for small data  in critical Besov spaces $\dot{B}^{3/p-1}_{p, \infty}$  with $3<p\le 6$. Chemin \cite{Che99} obtained global solutions in 3D for small data  in critical Besov spaces $\dot{B}^{3/p-1}_{p, q}$ for all $p<\infty, \ q\le \infty$.  Koch and Tataru  \cite{KoTa01} studied local solutions for initial data in $\mathit{vmo}^{-1}$ and global solutions for small initial data in $BMO^{-1}=\dot F^{-1}_{\infty,2}$. Iwabuchi \cite{Iw10} considered the local well posedness of NS for initial data in modulation spaces $M^{s}_{p,q}$, especially in $M^{-1}_{p,1}$ with $1\le p\le \infty$ and the global well-posedness with $p\leq d$.

On the other hand, NS is ill-posed in all critical Besov spaces $\dot{B}^{-1}_{\infty,q}$ with $1\leq q\leq \infty$. Indeed, Bourgain and Pavlovic \cite{BoPa08} obtained that NS is ill-posedness in $\dot{B}^{-1}_{\infty,\infty}$, Germain \cite{Ge08} and Yoneda \cite{Yo10} showed that NS is ill-posedness in $\dot{B}^{-1}_{\infty,q}$ for $q>2$. Finally, the ill-posedness of NS in all $\dot{B}^{-1}_{\infty,q}$, $1\leq q <\infty$ was obtained in \cite{Wa15}, where the ill-posedness means that the solution map $u_0\to u(t)$ is discontinuous at origin $t=0$.

Foias and Temam \cite{FoTe89} first proved spatial analyticity for the periodical solutions in Sobolev spaces. The analyticity of solutions in $L^p$ for NS was shown by Gruji\v{c} and Kukavica \cite{GrKu98}, and Lemari\'{e}-Rieusset \cite{LeRi00} gave a different approach based on multilinear singular integrals. Using iterative derivative estimates, the analyticity of NS for small initial data in $BMO^{-1}$ was obtained in Germain, Pavlovic and Staffilani \cite{GePaSt07} (see also Dong and Li \cite{DoLi09}, Miura and Sawada \cite{MiSa06} on the iterative derivative techniques). Bae, Biswas and  Tadmor \cite{BaBiTa12} obtained the analyticity of the solutions of NS in 3D for sufficiently small initial data in critical Besov spaces $\dot{B}^{3/p-1}_{p, q}$ with $1< p,q<\infty$.

Up to now, the well-posedness results of NS is known only for initial data in critical or subcritical function spaces. In super-critical function spaces, we have no well-posedness result on NS. For any $s<0$, modulation space $E^s_{p,q}$ is a super-critical space of NS in the sense that for any $\tilde{s}<-1+d/p$, $B^{\tilde{s}}_{p,q}$ is a super-critical space and $E^s_{p,q}$ is rougher than $B^{\tilde{s}}_{p,q}$, i.e., $B^{\tilde{s}}_{p,q} \subset E^s_{p,q}$. We can show that NS has a unique global solution if the initial data belong to $E^s_{r,1}$ ($s<0$) and their Fourier transforms are supported in one octant. Now we state our main results.

 \begin{thm}\label{NSthm2}
Let $d\ge 2$, $s<0$, $2\leq r \leq d$.  Assume that $u_0 \in   E^{s}_{r, 1} $ and ${\rm supp}\,\widehat{u}_0 \subset \mathbb{R}^d_I:=\{\xi:  \xi_i \geq 0, \ i=1,...,d\}$.   Then there exists $s_0 \leq s$  such that \eqref{NS3} has a unique mild solution $u \in  \widetilde{L}^{\gamma}(0,\infty; E^{s_0}_{p, 1}) \cap C (0,\infty; E^{s_0}_{r, 1}) $, where
$$
(\gamma, p) =
\left\{
\begin{array}{ll}
(d+2, d+2), & r=d\\
(2, \rho), &  2\leq r<d
\end{array}
\right.
$$
for some sufficiently large $\rho :=\rho(r)$  and
\begin{align}
      \|u\|_{\widetilde{L}^\gamma(I, E^s_{p,1})} =  \left\|\{2^{s|k|} \|\Box_k u\|_{L^\gamma_{t\in I}L^p_x} \}_{k\in \mathbb{Z}^d}\right\|_{\ell^1}.  \label{NSspace1}
\end{align}
Moreover, there exists a small constant $c>0$ such that the solution is real analytic when $t> s^2/c^2$.
\end{thm}

If $u_0$ is sufficiently small in $  E^{s}_{r, 1} $, we can take $s_0=s$ in Theorem \ref{NSthm2}.
Obviously, $L^2 \subset E^s_{2,1}$ for any $s<0$. So,  Theorem \ref{NSthm2} covers the $L^2$ initial data $u_0$ with ${\rm supp}\,\widehat{u}_0 \subset \mathbb{R}^d_I$.  Moreover, the initial condition $u_0 \in   E^{s}_{r, 1} $ can be replaced by
$$
u_0 \in E^{s_+}_{2,2}, \ \ \|u_0\|_{E^a_{2,2}} \sim \|2^{a|\xi|} \widehat{u}_0\|_2.
$$
Noticing that the partial derivatives of the Dirac measure $D^\alpha \delta \in E^s_{r,1}$ for any $s<0$, we see that for the initial data $u_0$ satisfying $\widehat{u}_0=    \overrightarrow{P}(\xi) \chi_{\mathbb{R}^d_I}$ or $\widehat{u}_0=   \overrightarrow{P}(\xi) 2^{-s|\xi|/2} \chi_{\mathbb{R}^d_I}$,  where   $\overrightarrow{P}(\xi)$ is a polynomial and $\overrightarrow{P}(\xi)\cdot \xi=0$, the solution of NS is globally existing, unique and real analytic for any $t>t_s$, where $t_s=s^2/c^2$.

Unfortunately, ${\rm supp}\,\widehat{u}_0 \subset \mathbb{R}^d_I$ implies that $u_0$ is not real valued. One may further ask if the condition ${\rm supp}\,\widehat{u}_0 \subset \mathbb{R}^d_I$ is necessary.  Using the ideas in showing the ill-posedness of NS in  $B^{-1}_{\infty, q}$, we can easily show that   ${\rm supp}\,\widehat{u}_0 \subset \mathbb{R}^d_I$ cannot be removed for the global well-posedness in $E^s_{r,1}$  (cf. \cite{Wa15}).

 As a further topic we study the global well-posedness of NS in modulation spaces $M^{-1}_{r,1}$, Iwabuchi \cite{Iw10} already obtained the global well-posedness of NS with small data in $M^{-1}_{r,1}$, $r\leq d$.  We can generalize his result to the case $d<r< \infty$ and consider the initial data in $\dot M^{-1}_{r, 1}$.

\begin{thm}\label{NSthm1}
Let $d\ge 2$, $d\leq r< p \leq 2r$.  Assume that $u_0 \in  \dot M^{-1}_{r, 1}$ is sufficiently small.  Then  \eqref{NS3} has a unique mild solution $u \in  \widetilde{L}^{2}(0,\infty; M^0_{p, 1}) \cap L^\infty(0,\infty; \dot M^{-1}_{r, 1}) $, where
\begin{align}
      \|u\|_{\widetilde{L}^\gamma(I, M^s_{p,1})} =  \left\|\{\langle k\rangle^s \|\Box_k u\|_{L^\gamma_{t\in I}L^p_x} \}_{k\in \mathbb{Z}^d}\right\|_{\ell^1}.  \label{2NSspace1}
\end{align}

\end{thm}

Throughout this paper, we will use the following notations. $C\ge 1, \ c\le 1$ will denote constants which can be different at different places, we will use $A\lesssim B$ to denote   $A\leqslant CB$; $A\sim B$ means that $A\lesssim B$ and $B\lesssim A$ and sometimes we write $A\lesssim_\lambda B$ to emphasize that the constant in $A\leq  C(\lambda) B$ depends on a parameter $\lambda$. We write $a\vee b= \max (a,b)$ and $a\wedge b = \min(a,b)$, $s_+= s+ \varepsilon$ with $0<\varepsilon\ll 1$. $\chi_E$ denotes the characteristic function on $E\subset \mathbb{R}^d$.   We denote  $|x|_\infty=\max_{i=1,...,d} |x_i|$, $|x|=|x_1|+...+|x_d|$ for any $x\in \mathbb{R}^d$. $L^p=L^p(\mathbb{R}^d)$ ($\ell^p$) stands for the (sequence) Lebesgue space for which the norm is written as $\|\cdot\|_p$.  Let $X$ be a Banach space. For any $I\subset \mathbb{R}^+=[0,\infty)$,   we denote
$$
\|u\|_{ L^{\gamma}(I; \ X)} = \left(\int_I \|u(t, \cdot)\|^{\gamma}_{X} dt \right)^{1/{\gamma}}
$$
 for $1\leq \gamma <\infty$ and with usual modifications
for $\gamma =\infty$. In particular,  if $X=L^p$,  we will write
$\|u\|_{L_{t\in I}^{\gamma}L_x^p}=\|u\|_{L^{\gamma}(I;L^p)}$ and $\|u\|_{L_{t }^{\gamma}L_x^p}=\|u\|_{L^{\gamma}(0,\infty;L^p)}$.

The paper is organized as follows. In Section 2 we review some known results on Besov and Trieble spaces, and $L^p$-multiplier, which will be frequently used in this paper. In Section 3 we mainly consider the relations between Gelfand-Shilov's space and the modulation space with exponential weights. We introduce a hybrid Riesz potential-modulation space in Section 4, which is useful in the study of NS. In Section 5 we will show an equivalent norm of modulation spaces with negative exponential weight by using the frequency uniform decomposition operators, their relations to Besov spaces are established.  Section 6 is devoted to the study of the semigroup $e^{t\Delta}$ and some basic estimates are obtained.  Theorem \ref{NSthm1} will be proven in Section 7. By considering a dilation  property of $E^s_{p,q}$ in Section 8, one can scale the large data of NS in the space $E^s_{p,q}$ $(s<0)$ to sufficiently small data. Theorem 1.1 will be shown in Sections 9 and 10. In the Appendix we give the proofs of some results which are known for the experts engaged in the study of modulation spaces but may be unknown for the specialists in the other subjects.

\section{Preliminary on Besov and Triebel spaces, $L^p$-multiplier}

Recall the definition of dyadic decomposition in
Littlewood-Paley theory \cite{Tr83}. Let $\psi: \mathbb{R}^d \rightarrow
[0, 1]$ be a smooth cut-off function which equals $1$ on the unit ball and equals $0$ outside the ball $\{\xi
\in \mathbb{R}^d; |\xi|\leq 2\}$.
Write $\varphi(\xi):=\psi(\xi)-\psi(2\xi)$ and
$\varphi_j(\xi)=\varphi(2^{-j}\xi)$.
$\Delta_j:=\mathscr{F}^{-1}\varphi_j \mathscr{F}, j\in \mathbb{Z}$
are said to be the dyadic decomposition operators which satisfy the
operator identity: $I=\sum_{j=-\infty}^{+\infty}\Delta_j$.  We write $\dot {\mathscr{S}}= \{\varphi\in \mathscr{S}: D^\alpha \widehat{\varphi}(0)=0, \ \forall \ \alpha \in \mathbb{Z}^d\}$ and  $\dot{\mathscr{S}'}$ denotes its dual space.
 The norms in homogeneous Besov spaces are defined as the subspace of $\dot {\mathscr{S}'}$ for which the norm is defined by
\begin{align}
\|f\|_{\dot{B}^s_{p, q}}= \|\{2^{sj} \Delta_j f\}_{j\in \mathbb{Z}}\|_{\ell^q(L^p)}=
 \left(\sum_{j=-\infty}^{+\infty}2^{jsq}\|\Delta_j f\|^q_{p}\right)^{1/q}
\label{Besov}
\end{align}
with usual modification if $q=\infty.$ We can define the Besov spaces $B^s_{p,q} = L^p\cap \dot B^s_{p,q}$ for $s>0$ and $B^s_{p,q} = L^p + \dot B^s_{p,q}$ for $s<0$.  Similarly, one can define the Triebel spaces (cf. \cite{Tr83}):
\begin{align}
\|f\|_{\dot{F}^s_{p, q}}= \|\{2^{sj} \Delta_j f\}_{j\in \mathbb{Z}}\|_{L^p(\ell^q)}
\label{Triebel}
\end{align}
We have $\dot{F}^s_{p, p\wedge q} \subset \dot{F}^s_{p, q} \subset \dot{B}^s_{p, p\vee q} $, $\dot{F}^{s_0}_{p_0, q} \subset \dot{F}^{s_1}_{p_1, r}$ and $\dot{F}^{s_0}_{p_0, \infty} \subset   \dot{B}^{s_1}_{p_1, p_0}$ if $s_0-d/p_0=s_1-d/p_1$, $s_0>s_1$. We have the following Gagliardo-Nirenberg inequality (see \cite{HaMoOzWa11}): Let $0<\theta <1$ and $s=(1-\theta)s_0+\theta s_1$. Then $\dot B^{s_0}_{p,\infty} \cap \dot B^{s_1}_{p,\infty} \subset \dot B^{s}_{p,1}$ and
\begin{align}
\|f\|_{\dot B^{s}_{p,1}} \lesssim  \|f\|^{1-\theta}_{\dot B^{s_0}_{p,\infty}}  \|f\|^{\theta}_{\dot B^{s_1}_{p,\infty}}.  \label{GNGN}
\end{align}
Let recall that Riesz potential spaces $\dot H^s_p := (-\Delta)^{-s/2} L^p$, $\dot H^s := H^s_2$ and we have
$$
\dot H^{s}_p =\dot F^{s}_{p  ,2}.
$$

As the end of this section, we recall a criterion of the multiplier on $L^p$.  We denote by $M_p$ the multiplier space on $L^p$, i.e.,
$$
\|m\|_{M_p} = \sup_{\|f\|_p=1} \|\mathscr{F}^{-1} m \mathscr{F} f\|_p.
$$
The following Bernstein's multiplier estimate is well-known: for $L>d/2$, $L\in \mathbb{N}$, $\theta= d/2L$ (cf. \cite{BL76, Ho60, WaHuHaGu11}),
\begin{align} \label{multiplier}
\|m\|_{M_p} \lesssim \|\mathscr{F}^{-1} m\|_1 \lesssim \|m\|^{1-\theta}_2 \|m\|^{\theta}_{\dot H^L}.
\end{align}

\section{Modulation spaces with tempered ultra-distributions}

The short-time Fourier transform (STFT) of a function $f$ with respect to a window function $g \in
\mathscr{S}$ is defined as follows (see \cite{Fei83, Groch01}):
\begin{align} \label{STFT}
V_g f(x, \xi) = \int_{\mathbb{R}^d} e^{-{\rm i} t \cdot \xi}
\overline{g(t-x)} f(t)dt  = \langle f, \ M_\xi T_x g\rangle,
\end{align}
where  $T_x:  g \to g(\cdot-x)$ and $  M_\xi : g  \to e^{\mathrm{i}\cdot  \, \xi} g$ denote the translation and modulation operators, respectively. The STFT is closely related to the wave packet transform of C\'{o}rdoba and Fefferman \cite{CoFe78, KaKoIt14} and the Wigner transform \cite{Groch01, Te06}. It is  a basic tool in time frequency analysis theory. It is known that (cf. \cite{Groch01})
\begin{align} \label{STFTprop}
V_g f(x, \xi) = e^{-i x\xi} V_{\hat{g}} \widehat{f} (\xi, -x) = e^{-i x\xi} (\mathscr{F}^{-1} T_{\xi} \overline{\hat{g}} \widehat{f}) (x).
\end{align}
Let $1\leq p,q\leq \infty$, $s\in \mathbb{R}$. The STFT equipped with the $L^q_\xi L^p_x$ norm  generates modulation spaces $M^s_{p,q}$  for which the norm is given in \cite{Fei83}:
\begin{align} \label{amodspace}
\|f\|^{\circ}_{M^s_{p,q}} =  \left\| \langle\xi\rangle ^{s} \|V_g f(x, \xi)\|_{L^p(\mathbb{R}^d_x)}\right\|_{L^q(\mathbb{R}^d_\xi)}
\end{align}
with the usual modifications if $p$ or $q$ is infinite. $M^s_{p,q}$ is defined as the
space of all tempered distributions $f\in \mathscr{S}' $ for
which $\|f\|^{\circ}_{M^s_{p,q}}$ is finite.   The notion of coorbit spaces were introduced in \cite{FeGr88} and as an example of coorbit spaces, modulation spaces were also generalized to cover some exponentially weighted
cases, i.e., $\langle\xi\rangle ^{s}$  is replaced by $e^{w (x,\xi)}$ in \eqref{amodspace}.  In this paper we mainly interested in the modulation spaces with negative exponential weights, which contain a class of distributions which cannot be handled by the Bj\"orck's tempered ultradistribution space.

\subsection{Tempered ultradistributions}

  Generalized distribution functions have a rich physical background, such as white noises come from physical phenomena which were naturally formulated as distribution function theory in stochastic processes (cf. \cite{De73, Ja76}). In current paper we mainly interested in the distribution functions for which their Fourier transforms have a definite meaning.  Since the convolution of a tempered distributions in $\mathscr{S}'$ and a Schwartz function has at most polynomial growth in $\mathbb{R}^d$, it seems interesting to make a generalization to the cases with exponential growth. A class of tempered ultradistributions were introduced by Bj\"orck \cite{Bj66} (see also Obiedat \cite{Ob06}).  We denote by $\mathfrak{M} $ the collection of all real-valued functions $\omega$ such that $\omega(\cdot) = \sigma(|\cdot|)$, where $\sigma (t)$ is an increasing continuous concave function on $[0,\infty)$ with $\sigma(0)=0$ and
\begin{align}
\int_{\mathbb{R}_+} \frac{\sigma(t)}{1+t^2} dt < \infty, \ \ \sigma (t) \geq c+ d \ln (1+t), \  \  \forall t\geq 0, \nonumber
\end{align}
where $c\in \mathbb{R}$  and $d>0$ are suitable constants. Let $\omega\in \mathfrak{M}$ and we denote by $\mathscr{S}_\omega $ the set of infinitely differentiable complex-valued functions $f$ satisfying
\begin{align}
p_{\lambda} (f) = \sup_{x\in \mathbb{R}^d} e^{\lambda \omega(x)} | f(x)|< \infty, \ \ \   q_{\lambda} (f) = \sup_{\xi \in \mathbb{R}^d} e^{\lambda \omega(\xi)} |\widehat{f}(\xi)|<\infty. \label{seminorm}
\end{align}
$\mathscr{S}_\omega$ equipped with the system of semi-norms $\{p_{\lambda}, \ q_{\lambda}\}_{\lambda>0}$  is a complete locally convex topological space.  We denote by $\mathscr{S}'_\omega$ the dual space of $\mathscr{S}_\omega$, which is said to be a tempered ultradistribution space up to sub-exponential growth, see \cite{Bj66, GrZi04}.

If $\omega(x) \sim \ln (1+|x|)$, then $\mathscr{S}_{\omega}= \mathscr{S}$. If $\omega(x) =|x|^\theta, \ \theta \in (0,1)$, then $\mathscr{S}_{\omega} \subset   \mathscr{S}$ and $\mathscr{S}'_\omega \supset \mathscr{S}' $ contains a class of distributions out of $\mathscr{S}'$.
However, if  $\omega(x) =|x|$, we have $\omega \notin \mathfrak{M}$.  \eqref{seminorm} becomes
\begin{align}
p_{\lambda,1} (f) = \sup_{x\in \mathbb{R}^d} e^{\lambda |x|} | f(x)|< \infty, \ \ \   q_{\lambda,1} (f) = \sup_{\xi \in \mathbb{R}^d} e^{\lambda |\xi|} |\widehat{f}(\xi)|<\infty. \label{seminorm1}
\end{align}
we denote by $\mathscr{S}_{1} $ the set of infinitely differentiable complex-valued functions $f$ satisfying $p_{\lambda,1} (f)<\infty$ and $q_{\lambda,1} (f) <\infty$.
$\mathscr{S}_{1} $ equipped with the system of semi-norms $\{p_{\lambda,1}, \ q_{\lambda,1}\}_{\lambda>0}$  is a complete locally convex topological space, which is said to be the Gelfand-Shilov space of Beurling type, cf. \cite{GeSh68, Si67, Pi88}.

\subsection{Gelfand-Shilov space via the STFT }

The modulation spaces with exponential weight  $e^{\lambda(\omega(x)+\omega(\xi))}$ with $\omega \in \mathfrak{M}$ in \eqref{amodspace} have been studied in \cite{GrZi04}.  Teofanov \cite{Te06} studied the case $\omega(\cdot)= |\cdot|^{\gamma}$ and the corresponding modulation spaces with $0<\gamma<1$.
 For our purpose we only consider the case $\omega(x) =|x|$. For convenience, we denote
$$
\|F\|_{L^{m}_{p,q}} = \left(\int_{\mathbb{R}^{d}} \left(\int_{\mathbb{R}^{d}} |F(x,\xi)m  (x,\xi) |^p dx\right)^{q/p}d\xi \right)^{1/q}
$$
with usual modification in the case $p,q=\infty $  and
\begin{align} \label{modspace}
\|f\|_{M^{m}_{p,q}} : = \|V_g f \|_{L^{m}_{p,q}},
\end{align}
where $g:= \pi^{-d/2} e^{-|x|^2/2}$ and $V_g$ is the STFT in \eqref{STFT}. Denote
$$
m_\lambda (x,\xi)= e^{\lambda (|x| + |\xi|)}.
$$
We point that $M^{m_\lambda}_{1,1}$ is an example of the coorbit spaces in \cite{FeGr88}.
It was shown in \cite{GrZi04} (see also \cite{ChChKi96}) that $f\in \mathscr{S}_1$ is equivalent to  $f\in M^{m_\lambda}_{\infty, \infty}$ for all $\lambda>0$. Some recent studies on the
 pseudo-differential operators on the Gelfand-Shilov spaces related to the STFT can be founded in \cite{Te06a, To12, CaTo17, KoPi16} and references therein.

 Following Proposition 3.12 as in \cite{GrZi04} (see also \cite{Ch99}), one sees that  the the system of semi-norms on $ \mathscr{S}_{1} $ can be replaced by $ \{ \|\cdot\|_{M^{m_\lambda}_{\infty,\infty}} \}_{\lambda>0}$. Following this fact, we can show that

\begin{prop} \label{GelfandShilov}
We denote $\mathscr{S}_{m} =\bigcap_{\lambda>0} M^{m_\lambda}_{1,1}$, which equipped with the system of norms $\{\|\cdot\|_{M^{m_\lambda}_{1,1}}\}_{\lambda>0}$  generates a complete locally convex topological space. We have $\mathscr{S}_m=\mathscr{S}_1$ with equivalent topologies.
\end{prop}
{\bf Proof.} Noticing that $M^{m_\lambda}_{1,1} \subset M^{m_\lambda}_{\infty,\infty} \subset  M^{m_{\lambda_+}}_{1,1}$ and $\mathscr{S}_{m} =\bigcap_{\lambda>0} M^{m_\lambda}_{\infty, \infty}$ (cf. \cite{GrZi01}), we have the result, as desired.  See also Appendix for the details of the proof by following \cite{Ch99}. $\hfill\Box$\\

We easily see that $\mathscr{S}_{m} =\mathscr{S}_{1}$ contains at least Gaussian and the finite linear combinations of Hermite functions. Since $\mathscr{S}_{m}$ is invariant under the translation and modulation, it contains also $M_\xi T_x e^{-|t|^2/2}$ and their linear combinations.  One may further ask if we can choose $m_\lambda = e^{\lambda (|x|^{1+\varepsilon} + |\xi|^{1+\varepsilon})}$, $\varepsilon>0$, in the case $\varepsilon<1$ it is possible, however, if $\varepsilon=1, \ \lambda>1/2$, then $ M^{m_\lambda}_{1,1}$ contains only the function $0$, cf. \cite{BoDeJa03, GrZi01, GrZi04}. One may further ask if we can replace $ M^{m_\lambda}_{1,1}$ by $ M^{m_\lambda}_{p,q}$, the answer is yes, however, noticing that $ M^{m_\lambda}_{1,1}\subset M^{m_\lambda}_{p,q}$,  the space $\mathscr{S}'_{m}$ seems more natural if we choose $\{M^{m_\lambda}_{1,1}\}_{\lambda>0}$ to generate  $\mathscr{S}_{m}$ than the other spaces $\{M^{m_\lambda}_{p,q}\}_{\lambda>0}$, $p$ or $q\neq 1$.

We denote by $\mathscr{S}'_{m}= \mathscr{S}'_{1}$ the dual space of $\mathscr{S}_{1}$.  By Proposition  \ref{GelfandShilov}  we have

\begin{cor}
Let $\omega \in \mathfrak{M}$. Then $\mathscr{S}'  \subset   \mathscr{S}'_{\omega} \subset \mathscr{S}'_m  $ are continuous embeddings.
\end{cor}

 Using a standard argument, we see that

\begin{lem} \label{Sprime} We have
$$
\mathscr{S}'_m = \bigcup_{\lambda>0} (M^{m_{\lambda}}_{1,1})^*,
$$
where  $(M^{m_{\lambda}}_{1,1})^*$ denotes the dual space of $ M^{m_{\lambda}}_{1,1} $. Namely, $f\in \mathscr{S}'_m$ if and only if there exists $\lambda>0$ and $C_\lambda >0$ such that
\begin{align}
|\langle f, \, \varphi\rangle| \leq C_\lambda \|\varphi\|_{M^{m_{\lambda}}_{1,1}}, \ \ \forall \, \varphi \in  \mathscr{S}_m. \label{continuousfunct}
\end{align}
\end{lem}
{\bf Proof.} It is easy to see that  $\mathscr{S}_m$ can be generated by the sequence of norms $\{\|\cdot\|_{M^{m_n}_{1,1}}\}^\infty_{n=1}$. We show that  $f\in \mathscr{S}'_m $  implies that \eqref{continuousfunct} holds for all $\lambda= n\in \mathbb{N}$.   If  not, then there exists $f\in \mathscr{S}'_m$ and $\varphi_n \in \mathscr{S}_m $ varifying
\begin{align}
|\langle f, \, \varphi_n \rangle| \geq n  \|\varphi_n \|_{M^{m_n}_{1,1}} . \label{continuousfunct1}
\end{align}
Denote $\psi_n =\varphi_n/(n \|\varphi_n \|_{M^{m_n}_{1,1}})$. Then for any $\varepsilon>0$, let $1/n_0 < \varepsilon$. By $\|\cdot\|_{M^{m_n}_{1,1}} \leq \|\cdot\|_{M^{m_{n'}}_{1,1}}$ for $n\leq n'$, we have
$\psi_n \in \{\psi \in \mathscr{S}_m:  \|\psi\|_{M^{m_{n_0}}_{1,1}} < \varepsilon\}$ if $n>n_0$.  Hence $\psi_n \to 0$ in $\mathscr{S}_m$, but
\begin{align}
|\langle f, \, \psi_n \rangle| \geq 1 . \label{continuousfunct2}
\end{align}
This contradicts with $f\in \mathscr{S}'_m$. By the density of $\mathscr{S}_m$ in $M^{m_\lambda}_{1,1}$ (see \cite{Groch01}, or Appendix), we have that $f$ satisfying \eqref{continuousfunct} is equivalent to  $f\in (M^{m_{\lambda}}_{1,1})^*$.      $\hfill\Box$ \\

Further, we can describe the dual space of $ M^{m_{\lambda}}_{1,1} $ by the STFT.
We claim that $V_g f $ is pointwise  defined for any $(x,\omega) \in \mathbb{R}^{2d}$. Indeed, by duality we have
$$
|V_g f (x,\omega) | =|\langle f, M_\omega T_x g\rangle| \leq \|f\|_{(M^{m_{\lambda}}_{1,1})^*} \|M_\omega T_x g\|_{M^{m_{\lambda}}_{1,1}}.
$$
Since $|V_g (M_{\omega_0} T_{x_0} g) (x,\omega)|= e^{-(|x-x_0|^2 + |\omega-\omega_0|^2)/4}$, we see that
$$
\|M_\omega T_x g\|_{M^{m_{\lambda}}_{1,1}} \lesssim_{\lambda}  e^{ \lambda(|x|+|\omega|)}.
$$
Hence,
$$
|V_g f(x,\omega)|  \lesssim_\lambda e^{ \lambda(|x|+|\omega|)}  \|f\|_{(M^{m_{\lambda}}_{1,1})^*} .
$$
It follows that $V_g f \in L^{1/m_\lambda}_{\infty,\infty}$ if $f\in (M^{m_{\lambda}}_{1,1})^*$.  Moreover, we can further describe  $(M^{m_{\lambda}}_{1,1})^*$ as follows.

\begin{lem}  {\rm (\cite{Groch01})}  \label{STFTcon}
Let $f\in \mathscr{S}'_m$. Then we have $V_g f $ is a continuous function on $\mathbb{R}^{2d}$.
\end{lem}

Recall that $M^{1/m_{\lambda}}_{\infty,\infty}$ is defined as the space of all $f\in \mathscr{S}'_1$ such that
$$
\|f\|_{M^{1/m_{\lambda}}_{\infty,\infty}} := \|e^{-\lambda(|x|+|\xi|)} V_g f\|_{L^{\infty}(\mathbb{R}^{2d})} < \infty.
$$

\begin{prop} [Duality] \label{duality}  We have $(M^{m_{\lambda}}_{1,1})^* =M^{1/m_{\lambda}}_{\infty,\infty}$ under the duality
\begin{align} \label{Dualint}
\langle f, \varphi\rangle =\int_{\mathbb{R}^{2d}} \overline{V_g f (x,\omega)} V_g \varphi (x,\omega) dx d\omega
\end{align}
for any $f\in M^{1/m_{\lambda}}_{\infty,\infty}$ and $\varphi \in M^{m_{\lambda}}_{1,1}$.
\end{prop}
{\bf Proof.} Let us follow \cite{Groch01}, Theorem 11.3.6. For completeness we give the details of the proof.  By H\"older's inequality, we see that any  $f\in M^{1/m_{\lambda}}_{\infty,\infty}$ generates a continuous functional on $M^{m_{\lambda}}_{1,1}$ via \eqref{Dualint}.

Next, $M^{m_{\lambda}}_{1,1}$ is isometric to a subspace of $L^{m_\lambda}_{1,1}$ and
$$
V_g : M^{m_\lambda}_{1,1} \to V= \{V_g \varphi:  \  \varphi\in M^{m_\lambda}_{1,1}\} \subset L^{m_\lambda}_{1,1}
$$
is an isometric mapping. So, each continuous functional $f \in (M^{m_\lambda}_{1,1})^*$ induces a bounded linear functional $\widetilde{f}$ on $V$ by letting
$$
\langle \widetilde{f}, V_g \varphi \rangle =  \langle  f ,   \varphi \rangle.
$$
Then $\widetilde{f}$ can be extended to continuous functional on $L^{m_\lambda}_{1,1}$ whose norm will be preserved (its extension is still written as $\widetilde{f}$). By the duality of $(L^{m_\lambda}_{1,1})^*= L^{1/m_\lambda}_{\infty,\infty}$,  there exists  $H\in L^{1/m_\lambda}_{\infty,\infty}$ such that,
$$
\langle \widetilde{f}, \psi \rangle = \int_{\mathbb{R}^{2d}} \overline{H(x,\omega)} \psi(x,\omega) dx d\omega, \ \ \forall \ \psi \in  L^{m_\lambda}_{1,1}.
$$
It follows that
$$
\langle  f ,   \varphi \rangle =  \int_{\mathbb{R}^{2d}} \overline{H(x,\omega)} V_g\varphi(x,\omega) dx d\omega, \ \  \forall \  \varphi\in M^{m_\lambda}_{1,1}. $$
Now we introduce
$$
V^*_g H = \int_{\mathbb{R}^{2d}} H(x,\omega)  M_\omega T_x g dx d\omega.
$$
We can show that $V^*_g H  \in M^{1/m_\lambda}_{\infty,\infty}$. In fact, we easily see that
$$
V_g(V^*_g H) \leq |H|* |V_g g|.
$$
It follows from $1/m_{\lambda}(y) \leq m_\lambda (y-x)/ m_{\lambda}(x)$ and Young's inequality that
$$
\|V^*_g H\|_{M^{1/m_\lambda}_{\infty,\infty}} \leq \|H\|_{L^{1/m_\lambda}_{\infty,\infty}} \|V_g g\|_{L^{m_\lambda}_{1,1}}.
$$
So, it follows from Fubini's theorem that for any $\varphi \in M^{m_\lambda}_{1,1} $,
$$
\langle V^*_g H, \varphi \rangle =  \int_{\mathbb{R}^{2d}} \overline{H(x,\omega)}  V_g \varphi (x,\omega)  dx d\omega = \langle f, \varphi \rangle .
$$
It follows that $f= V^*_g H$.
We have the result, as desired. $\hfill\Box$

\begin{rem} \rm Similarly, we have $(M^{m_{\lambda}}_{p,q})^* =M^{1/m_{\lambda}}_{p',q'}$ for all $1\leq p,q<\infty$.  The integral in \eqref{Dualint} can be understood as the limit of Riemann sums, since both $V_g f$ and $V_g \varphi $  are continuous functions. More details related to this question can be found in \cite{FeWe06}.

\end{rem}

\begin{cor} {\rm (\cite{Groch01})}
Let $f\in  \mathscr{S}'_m$. Then  there exist $\lambda>0$ such that
$$
|V_g f (x,\omega)| \lesssim e^{\lambda(|x|+|\omega|)} .
$$
\end{cor}

\begin{lem}  \label{PropSm}
 Let $\alpha\in \mathbb{Z}^d$ and $X^\alpha f (t) = t^\alpha f(t)$. We have the following results:
\begin{itemize}
\item[\rm (i)]
  $ \mathscr{F}:  \mathscr{S}_m  \to  \mathscr{S}_m $ is an isomorphism.

\item[\rm (ii)]  $D^\alpha: \mathscr{S}_m \to \mathscr{S}_m$ is a continuous mapping.

\item[\rm (iii)]  $X^\alpha: \mathscr{S}_m \to \mathscr{S}_m$ is a continuous mapping.
  \end{itemize}
\end{lem}
{\bf Proof.} Noticing that   $|V_g f (x,\omega)|= |V_{\widehat{g}} \widehat{f} (\omega, -x)|$, we immediately have
$$
\|f\|_{M^{m_\lambda}_{1,1}} = \|\widehat{f}\|_{M^{m_\lambda}_{1,1}}.
$$
It follows that $\mathscr{F}: M^{m_\lambda}_{1,1} \to  M^{m_\lambda}_{1,1} $ is an isometric mapping, which implies that  $ \mathscr{F}:  \mathscr{S}_m  \to  \mathscr{S}_m $ is an isomorphism. To prove (ii), let us observe that for some polynomial $P(\cdot)$,
\begin{align} \label{Idineq}
e^{\lambda |t|}|D^\alpha f(t)| \lesssim  e^{  \lambda|t|} \int_{\mathbb{R}^{2d}} |V_g f |  |P(\omega)| e^{-|t-x|^2/4}  dxd\omega \lesssim \|f\|_{M^{m_\lambda}_{1,1}}.
\end{align}
Together with Proposition \ref{GelfandShilov},  we can show that the result of (ii) and (iii) and the details are omitted. $\hfill\Box$\\

By Lemma \ref{PropSm}, we can define the Fourier transform and multi-derivatives in $\mathscr{S}'_m$  exactly as in $\mathscr{S}'$. For any $f\in \mathscr{S}'_m, \  \varphi\in \mathscr{S}_m$, we define
$$
\langle \mathscr{F}f, \varphi \rangle :=  \langle f, \mathscr{F}\varphi \rangle, \ \ \ \  \langle D^\alpha f, \varphi \rangle :=  \langle f,  (-D)^\alpha\varphi \rangle,
$$
then we can make the Fourier transform and derivative operations in $\mathscr{S}'_m$.  Further, we can define the convolution operation in $\mathscr{S}_m$:
$$
f*g(t) =\int_{\mathbb{R}^{d}} f(t-y)g(y) dy.
$$

\begin{lem} [Convolution algebra] \label{Propconv}
  We have $*: L^{m_\lambda}_{1,1} \times  M^{m_\lambda}_{1,1} \to  M^{m_\lambda}_{1,1}$  and so,  $*: \mathscr{S}_m \times \mathscr{S}_m \to \mathscr{S}_m$ are continuous mappings. More precisely, we have
$$
\| \varphi*\psi \|_{M^{m_\lambda}_{1,1} }  \lesssim   \| \varphi\|_{M^{m_\lambda}_{1,1}}  \| \psi\|_{L^{m_\lambda}_{1,1}}   \lesssim_\lambda  \| \varphi\|_{M^{m_\lambda}_{1,1}}  \| \psi\|_{M^{m_\lambda}_{1,1}}.
$$
\end{lem}

{\bf Proof.}  We have for any $\varphi, \ \psi \in M^{m_\lambda}_{1,1}$,
$$
V_g (\varphi*\psi)(x, \omega) =\int_{\mathbb{R}^{d}} V_g \varphi (x-y, \omega) e^{{\rm i} y \omega} \psi(y) dy.
$$
Using $e^{\lambda|x|} \leq  e^{\lambda|x-y|} e^{\lambda|y|}$ we have
$$
\|V_g (\varphi*\psi)\|_{L^{m_\lambda}_{1,1} }  \leq \|V_g \varphi\|_{L^{m_\lambda}_{1,1}}  \int_{\mathbb{R}^{d}}  e^{\lambda|y|}  |\psi(y)| dy.
$$
Since
$
\psi = \int_{\mathbb{R}^{2d}} V_g \psi (x,\omega) M_\omega T_x g dx d\omega,
$
we have
\begin{align}
\int_{\mathbb{R}^{d}}  e^{\lambda|y|}  |\psi(y)| dy &  \lesssim   \int_{\mathbb{R}^{2d}} e^{\lambda|x|} |V_g \psi (x,\omega)| dxd\omega \int_{\mathbb{R}^{d}} e^{\lambda|x-y|} e^{-|x-y|^2/2} dy  \nonumber\\
& \lesssim_\lambda  \|V_g \psi\|_{L^{m_\lambda}_{1,1}}  =  \|  \psi\|_{M^{m_\lambda}_{1,1}}.  \nonumber
\end{align}
So, we have the result, as desired. $\hfill \Box$\\

Following the proof of Lemma \ref{Propconv}, in order to guarantee that $\varphi*\psi \in M^{m_\lambda}_{1,1}$, it suffices to assume that $\varphi\in M^{m_\lambda}_{1,1}, \ \psi \in L^1(\mathbb{R}^d, e^{\lambda|x|}dx)$.  By Lemma \ref{Propconv}, we can define the convolution $\varphi* f$ for any $\varphi\in \mathscr{S}_m$, $f\in \mathscr{S}'_m$:
$$
\langle \varphi * f, \ \psi\rangle :=  \langle f, \  \varphi (-\cdot) * \psi\rangle.
$$

Noticing that $\widehat{\varphi\psi} = \widehat{\varphi } * \widehat{\psi}$  holds in $\mathscr{S}_m$, by Lemmas \ref{PropSm} and \ref{Propconv}, we can also define $\varphi  f$ for any $\varphi\in \mathscr{S}_m$, $f\in \mathscr{S}'_m$:
$$
\varphi  f = \mathscr{F}^{-1} (\widehat{\varphi } * \widehat{f}).
$$

Recall that $\mathscr{F}\varphi \in \mathscr{S}_\omega$ if $\varphi \in C^\infty_0$. Unfortunately, this property does not hold in $ \mathscr{S}_m$; cf. \cite{GrZi04}.  One may think that we can directly use $M^{m_\lambda}_{1,1}\subset \mathscr{S}_\omega$ as a test function space and its dual $M^{1/m_\lambda}_{\infty,\infty}$ as a tempered ultra distribution space, since most useful distributions have been included by $M^{1/m_\lambda}_{\infty,\infty}$ and moreover, a norm structure is much simpler than the topological structure of $\mathscr{S}'_m$.  However,  $D^\alpha $ cannot be a bounded operator in $M^{m_\lambda}_{1,1}$ and we can only have $D^\alpha: M^{m_{\lambda+\varepsilon}}_{1,1} \to M^{m_\lambda}_{1,1}$ for any $\varepsilon>0$, this is why we consider the dual of $\mathscr{S}_1=\mathscr{S}_m$ as a tempered ultradistribution space.

\section{Hybrid Riesz potential-modulation spaces}

 Let $\Box_k$ be as in \eqref{FUD}. $\Box_k$ are almost orthogonal, i.e., $\Box_k =   \Box_k \circ \sum_{|\ell|_\infty \leq 1}\Box_{k+\ell}$. For convenience, we will write $\widetilde{\sigma}_k = \sum_{|\ell|_\infty \leq 5}\sigma_{k+\ell}$ and $\widetilde{\Box}_k =   \sum_{|\ell|_\infty \leq 1}\Box_{k+\ell}$.

 Let $s\in \mathbb{R}$, $1< p < \infty, \ 1\leq q\leq \infty$.   We denote by $ f\in \dot M^s_{p,q}$ that $f$ is in $\dot{\mathscr{S}}'({\Bbb R}^d)$ for which the following norm is finite:
\begin{align}
  \|f\|_{\dot M^s_{p,q}} :=  \left\|\{  (-\Delta)^{s/2} \Box_k f\}_{k\in \mathbb{Z}^d}  \right\|_{\ell^q(L^p)}  \label{defmod1aa}
\end{align}
and $\dot M^s_{p,q}$  is said to be the hybrid Riesz potential-modulation space.

\begin{prop}\label{Homomod1} Let $s\in \mathbb{R}$, $1< p <\infty, \ 1\leq q\leq \infty$. Then we have
$$
\|f\|_{\dot M^s_{p,q}} \sim  \left\|\left(\sum_{j\lesssim 1} 2^{2sj} |\triangle_j f |^2 \right)^{1/2}\right\|_p + \left(\sum_{k\neq 0} \langle k\rangle^{sq}\|\Box_k f\|^q_p\right)^{1/q}.
$$
\end{prop}
{\bf Proof.} Noticing that in the definition of $\dot M^s_{p,q}$,
$$
\|(-\Delta)^{s/2} \Box_k f\|_p \sim |k|^{s} \|\Box_k f\|_p \sim \langle k\rangle^{s} \|\Box_k f\|_p, \ \ k\neq 0,
$$
we have
$$
\left\|\{  (-\Delta)^{s/2} \Box_k f\}_{k\neq 0}  \right\|_{\ell^q(L^p)} \sim \left(\sum_{k\neq 0} \langle k\rangle^{sq}\|\Box_k f\|^q_p\right)^{1/q}.
$$
For $k=0$, one sees that
$
\|(-\Delta)^{s/2} \Box_0 f\|_p =\|\Box_0 f\|_{\dot H^s_p}.
$
From the equivalent norm on $\dot H^s_p$ it follows that
$$
\|f\|_{\dot H^s_p} \sim  \left\|\left(\sum_{j\in \mathbb{Z}} 2^{2sj} |\triangle_j f |^2 \right)^{1/2}\right\|_p,
$$
which implies the result, as desired. $\hfill\Box$

So, we see that the lower frequency part of $\dot M^s_{p,q}$ is $\dot H^s_p= \dot F^s_{p,2}$, the higher frequency part of $\dot M^s_{p,q}$ is $M^s_{p,q}$.

\section{Modulation spaces with negative exponential weight} \label{Emodspace}

\subsection{Equivalent and complete norms, dual spaces}

Let $1\leq p,q\leq \infty$, $s\in \mathbb{R}$ and $m_s(\xi)= 2^{s|\xi|}$. Let us define
\begin{align}
M^{m_s}_{p,q} =\{f\in \mathscr{S}'_m: \|f\|_{M^{m_s}_{p,q}}: = \|m_s(\xi)V_gf(x,\xi)\|_{L^{p,q}} <\infty \}. \label{modspaceexp0}
\end{align}
Since $\mathscr{S}'_m$ is the collection of all $M^{1/m_\lambda}_{\infty, \infty}$, $m_\lambda(z) =e^{\lambda |z|}$,  $M^{m_s}_{p,q}$ is reasonable as a subspace of $\mathscr{S}'_m$ for all $s\in \mathbb{R}$.

From the PDE point of view, the solutions are usually easier to control by the discrete STFT. So, we are looking for the frequency-localized version of $M^{m_s}_{p,q}$, $m_s= 2^{s|\cdot|}$.  For any $0< s< \infty$, $0<p,q \le \infty$, the
following space
\begin{align} \label{modspaceexp}
E^s_{p,q}({\Bbb R}^d) =\left \{ f\in {\mathscr S}'({\Bbb R}^d):
\; \|f\|_{E^s_{p,q}} := \left(\sum_{k \in {\Bbb Z}^d}
(2^{s |k|} \|\,\Box_k f\|_p)^q \right)^{1/q}<\infty
\right\}
\end{align}
can be regarded as a class of modulation spaces with an exponential regularity, cf. \cite{Wa06}.
Unfortunately, if $s<0$, we cannot treat $E^s_{p,q}$ as a subspace of $\mathscr{S}'$. In fact, let us consider a simple case $s<0$, $p=2$. Put
$$
\widehat{f} (\xi) = 2^{-s_0|\xi|},  \ s< s_0 <0.
$$
It is easy to see that $f\in E^s_{2,q}$ for all $q\in [1,\infty]$, however, $f\not\in \mathscr{S}'$. Of course, $E^s_{2,q}$ contains any order partial derivatives of Dirac measure $D^\alpha \delta$. So, it is obliged to treat $E^s_{p,q}$ as a subspace of $\mathscr{S}'_m$ in the case $s<0$. On the other hand, $\mathscr{S}_{m} $ cannot contain a class of important functions whose Fourier transforms have compact support sets; cf. \cite{GrZi04}, which leads to that one cannot directly define $\Box_k f= (\mathscr{F}^{-1}\sigma_k ) * f$ for all $f\in \mathscr{S}'_m$, since $(\mathscr{F}^{-1}\sigma ) * f$ is only well defined for $\mathscr{F}^{-1} \sigma \in \mathscr{S}_m$.

One may consider the following sub-exponential regularity space  $E^s_{p,q, \theta} \subset \mathscr{S}'_\omega$ with $\omega(\cdot)= |\cdot|^\theta$, $\theta \in (0,1)$ for which the norm is given by
$$
\|f\|_{E^s_{p,q, \theta}} := \left(\sum_{k \in {\Bbb Z}^d}
(2^{s |k|^{\theta}} \|\,\Box_k f\|_p)^q \right)^{1/q}.
$$
$E^s_{p,q, \theta}$ $(s<0)$ can be well defined, since $\Box_k f$ is meaningful in $\mathscr{S}'_\omega$.  However, it is impossible to have any algebra structure for $E^s_{p,q, \theta}$ if $s<0$, even if we restrict the Fourier transform of of $f\in E^s_{p,q, \theta}$ supported in one octant.  Our solution is to introduce the following
\begin{align}
E^s_{p,q} = \big\{ f\in \mathscr{S}'_m: & \ \exists \ f_k \in L^p (\mathbb{R}^d), \ {\rm supp}\, \widehat{f}_k \subset k+ [-3/2, 3/2]^d  \nonumber\\
& {\rm such \ that} \ f=\sum_{k\in \mathbb{Z}^d} f_k \ {\rm in} \ \mathscr{S}'_m  \ {\rm and} \ \|\{2^{s|k|}f_k\}\|_{\ell^q(L^p)} <\infty \big\}, \label{espq}\\
\|f\|_{E^s_{p,q}} = & \inf  \|\{2^{s|k|}f_k\}\|_{\ell^q(L^p)}, \label{normespq}
\end{align}
where the infimum is taken over all of the decompositions of $f= \sum_{k\in \mathbb{Z}^d} f_k \in E^s_{p,q}$.

\begin{rem} \label{boxkiswell} \rm
Since $L^p \subset \mathscr{S}' \subset \mathscr{S}'_m$, we see that for any $f\in E^s_{p,q}$,  $ \widehat{f} = \sum_{k\in \mathbb{Z}^d} \widehat{f}_k$, $ {f}_k \in L^p$ implies that $\widehat{f}_k \in \mathscr{S}'$, so ${\rm supp}\,\widehat{f}_k$ is meaningful. In general, we cannot define $\Box_k f$ for all $f\in \mathscr{S}'_m$.  However, $\Box_k f$ is meaningful for any $f\in E^s_{p,q}$. Indeed,  we have $\lim_{N\to \infty} \sum_{|k|_\infty \leq N} \widehat{f}_k =\widehat{f}$ in $\mathscr{S}'_m$. In view of  $\sigma_k \sum_{|n|_\infty \leq N} \widehat{f}_n = \sigma_k \sum_{|\ell|_\infty \leq 5} \widehat{f}_{k+\ell} $ in $\mathscr{S}'$ for $N$ sufficiently large, one sees that $\sigma_k \widehat{f}= \sigma_k \sum_{|\ell|_\infty \leq 5} \widehat{f}_{k+\ell}$. It follows that $\Box_k f$ is well-defined if $f\in E^s_{p,q}$.

If we consider the regularity weight as the power functions $\langle \xi\rangle^s$, we see that modulation spaces $M^s_{p,q}$ have equivalent norms as in \eqref{defmod1} and \eqref{modspace}.  So, we conjecture that for the regularity weight as the exponential functions $2^{s|\xi|}$, we have similar equivalent norm in modulation spaces. However, $\Box_k f$ may have no meaning if $f\in \mathscr{S}'_m$ in \eqref{modspaceexp} and so, we cannot directly use the norm as in \eqref{modspaceexp} if $s<0$. Fortunately, we can use \eqref{normespq} as a substitute and we have
\end{rem}

\begin{prop}\label{inclusion1} Let $s\in \mathbb{R}$, $m_s(\xi)= 2^{s|\xi|}$,  $1\leq p,q\leq \infty$. Then we have
$E^s_{p,q}  = M^{m_s}_{p,q}$ with equivalent norms.
\end{prop}
{\bf Proof.} First, we show that $ M^{m_s}_{p,q} \subset E^s_{p,q}$.
For any $\xi \in k+ [-1/2, 1/2]^d$, we have from Young's inequality that
$$
\|\Box_k f\|_p \leq  \|\mathscr{F}^{-1} (\sigma_k e^{|\xi-\, \cdot|^2/2})\|_1  \|\mathscr{F}^{-1}  e^{-|\xi- \, \cdot|^2/2} \widehat{f} \|_p.
$$
It follows that $\Box_k f$ is well defined if $f\in M^{m_s}_{p,q}$.  For any $L>d/2$, we have for any $\xi \in k+[-1/2, 1/2]^d$,
$$
 \|\mathscr{F}^{-1} (\sigma_k e^{|\xi- \, \cdot|^2/2})\|_1 \leq  \| \sigma_k e^{|\xi-\, \cdot|^2/2} \|_{H^L} \lesssim 1.
$$
Hence, for $g= e^{-|\cdot|^2/2}$, we have for any $\xi \in  k+ [-1/2, 1/2]^d$,
$$
\|\Box_k f\|_p \leq  \|V_g f (\cdot,\xi) \|_{L^p(\mathbb{R}^d)}.
$$
It follows that
\begin{align}
\|f\|^q_{E^s_{p,q}} & \lesssim  \sum_{k\in \mathbb{Z}^d} \int_{k+[-1/2, 1/2]^d } 2^{s|\xi|q} \|\Box_k f\|^q_p d\xi.  \nonumber\\
&  \lesssim  \sum_{k\in \mathbb{Z}^d} \int_{k+[-1/2, 1/2]^d } 2^{s|\xi|q}    \|V_g f (\cdot,\xi) \|^q_{L^p(\mathbb{R}^d)} d\xi =  \|f\|^q_{M^{m_{s}}_{p,q}}. \nonumber
\end{align}
Hence, we obtain that $ M^{m_{s}}_{p,q} \subset E^s_{p,q}$.

Next, we prove that $ E^s_{p,q} \subset  M^{m_{s}}_{p,q} $. It suffices to show that
$$
\|2^{s|\xi|} V_g f\|_{L^q_\xi L^p_x(\mathbb{R}^{2d})} \lesssim \|f\|_{E^s_{p,q}}.
$$
holds for all $f\in E^s_{p,q}$.  Let us observe that
\begin{align}
\|  V_g f\|_{L^{m_s}_{p,q}(\mathbb{R}^{2d})}  \sim  \left(\sum_{k\in \mathbb{Z}^d} 2^{s|k|q} \int_{k+[-1/2, 1/2]^d } \|V_g f (\cdot,\xi) \|^q_{L^p(\mathbb{R}^d)} d\xi  \right)^{1/q}. \label{equivnorm2}
\end{align}
By $|V_gf|= |\mathscr{F}^{-1}\overline{\hat{g}}(\xi-\cdot) \widehat{f}|$ and $\widehat{f}= \sum_{m}\widehat{f}_m $,  we see that for any $\xi \in k+[-1/2, 1/2]^d$,
\begin{align}
\|V_g f\|_{L^p(\mathbb{R}^{d})} & \lesssim  \sum_{m\in \mathbb{Z}^d}  \| \mathscr{F}^{-1}     g(\xi-\cdot) \widehat{f}_m\|_{L^p(\mathbb{R}^{d})} \nonumber\\
& \lesssim     \sum_{m\in \mathbb{Z}^d} \left\|\mathscr{F}^{-1}   \left(g(\xi-\cdot ) \sum_{|\ell|_\infty \leq 6} \sigma_{m+\ell} \right) \right\|_{1}  \|{f}_m\|_{L^p(\mathbb{R}^{d})}. \nonumber
\end{align}
Using Bernstein's estimate, we have for any $\xi \in k+[-1/2, 1/2]^d$,
\begin{align}
\left\|\mathscr{F}^{-1}   (g(\xi-\cdot) \sum_{|\ell|_\infty \leq 6} \sigma_{m+\ell}) \right\|_{1}     \lesssim e^{-c|k-m|^2}. \nonumber
\end{align}
Hence, for any $\xi \in  k+ [-1/2, 1/2]^d$,
\begin{align}
\|V_g f(\cdot, \xi)\|_{L^p(\mathbb{R}^{d})} & \lesssim    \sum_{m\in \mathbb{Z}^d}    e^{-c|k-m|^2}  \| {f}_m \|_{p}.
\end{align}
Inserting the above inequality into \eqref{equivnorm2} and then using Young's inequality, we immediately obtain $ E^s_{p,q} \subset  M^{m_{s}}_{p,q} $. $\hfill\Box$

\begin{cor}\label{ainclusion1} Let $s\in \mathbb{R}$, $1\leq p,q\leq \infty$. Then we have
$$
\mathscr{S}_m \subset E^s_{p,q} \subset \mathscr{S}'_m
$$
with continuous embeddings.
\end{cor}

\begin{prop}\label{equivlent norm} Let $s\leq 0$, $1\leq p,q\leq \infty$. Denote
\begin{align}
\|f\|^{\circ}_{E^s_{p,q}} = \|\{2^{s|k|}\Box_k f\}\|_{\ell^q(L^p)}. \label{equivnorm}
\end{align}
Then $\|\cdot\|^{\circ}_{E^s_{p,q}}$ is an equivalent norm on $E^s_{p,q}$  (defined as in \eqref{espq}).
\end{prop}

{\bf Proof.} If $f\in E^s_{p,q}$, then by Remark \eqref{boxkiswell}, $\Box_kf$ has definite meaning. Let $ \|f\|^{\circ}_{E^s_{p,q}} <\infty$.  Following the proof of the second inclusion in Proposition \ref{inclusion1}, we have $f=\sum_k \Box_k f$ in $\mathscr{S}'_m$. Taking $f_k=\Box_k f$, one immediately gets that $\|f\|_{E^s_{p,q}} \leq \|f\|^{\circ}_{E^s_{p,q}}$.

Conversely, if
$$
f=\sum_k   f_k, \ \ f_k \in L^p, \ \ {\rm supp}\,\widehat{f}_k \subset k+[-3/2,3/2]^d,
$$
then by Remark \eqref{boxkiswell}, we have $\Box_k f = \sum_{|\ell|_\infty \leq 5} \Box_k {f}_{k+\ell}$. It follows that
$$
\|\Box_k f\|_p \leq \sum_{|\ell|_\infty \leq 5} \| {f}_{k+\ell}\|_p.
$$
Hence, $\|f\|^{\circ}_{E^s_{p,q}} \lesssim  \|\{2^{s|k|}f_k\}\|_{\ell^q(L^p)}. $ $\hfill\Box$ \\

One can also use the dilation of lattice points $\{\alpha k\}_{k\in \mathbb{Z}^d}$ with $\alpha>0$ to generate the equivalent norm.

\begin{prop}\label{equivalent norm} Let $s\leq 0$, $1\leq p,q\leq \infty$, $\alpha>0$. Denote $\Box_{\alpha k} = \mathscr{F}^{-1} \sigma_{\alpha k} \mathscr{F}$, $\sigma_{\alpha k}= \sigma(\cdot- \alpha k)$ and
\begin{align}
\|f\|^{(\alpha)}_{E^s_{p,q}} = \|\{2^{s\alpha| k|}\Box_{\alpha k} f\}\|_{\ell^q(L^p)}. \label{alphaequivnorm}
\end{align}
Then $\|\cdot\|^{(\alpha)}_{E^s_{p,q}}$ is an equivalent norm on $E^s_{p,q}$.
\end{prop}
{\bf Proof.} Using the almost orthogonality of $\Box_k$ and the multiplier estimate \eqref{multiplier}, we can get the result, as desired. $\hfill\Box$

In the following we give an equivalent norm in the case $1<p,q<\infty$.  Denote
$$
\Box^{c}_k =  \mathscr{F}^{-1} \chi_{k+[0,1)^d} \mathscr{F}.
$$
For $s\in \mathbb{R}$, $1<p,q<\infty$, we write
$$
\|f\|^{c}_{E^s_{p,q}} = \|\{2^{s|k|} \Box^{c}_k f\}\|_{\ell^q(L^p)}.
$$

\begin{prop}\label{equivnorm1} Let $s\in \mathbb{R}$, $1< p,q< \infty$. Then $\|\cdot \|^{\circ}_{E^s_{p,q}}$ and  $\|\cdot \|^{c}_{E^s_{p,q}}$ are equivalent norms on $E^s_{p,q}$. So, $\|\cdot \|_{E^s_{p,q}} \sim \|\cdot \|^{\circ}_{E^s_{p,q}} \sim \|\cdot \|^{c}_{E^s_{p,q}}$ on $E^s_{p,q}$.
\end{prop}

{\bf Proof.} Noticing that $\sigma_k \in M_p$ and $\chi_{k+[0,1)^d} \in M_p$ for any $p\in (1,\infty)$, we have
$$
\|\Box^c_k f\|_p \lesssim \sum_{|\ell|_\infty \leq 5} \|\Box_{k+\ell} f\|_p,  \ \ \|\Box_k f\|_p \lesssim \sum_{|\ell|_\infty \leq 5} \|\Box^c_{k+\ell} f\|_p,
$$
which implies the result, as desired. $\hfill\Box$

\begin{prop}\label{dual} Let $s\in \mathbb{R}$, $1\leq p,q< \infty$. Then we have
\begin{align}
(E^s_{p,q})^*  = E^{-s}_{p'q'}. \label{dualspace}
\end{align}
\end{prop}
{\bf Proof.} By Proposition \ref{inclusion1}, $E^s_{p,q} = M^{m_s}_{p,q}$ with equivalent norms.  Following the proof of Proposition \ref{duality}, we have $(M^{m_s}_{p,q})^* = M^{1/m_s}_{p',q'}$, see \cite{Groch01}.  $\hfill\Box$\\

\begin{prop}\label{Banachspace} Let $s\in \mathbb{R}$, $1\leq p,q\leq \infty$.
Then $\|\cdot\|_{E^s_{p,q}}$ is a complete norm on $E^s_{p,q}$, i.e. $E^s_{p,q}$ is a Banach space.
\end{prop}

{Proof.}  By $E^s_{p,q}= M^{m_s}_{p,q}$  and  $M^{m_s}_{p,q}$ is a dual space of $M^{1/m_s}_{p',q'}$ for $1<p,q<\infty$, it follows that $E^s_{p,q}$ is complete for $1<p,q<\infty$.  A direct proof in the cases that $p$ or $q$ equals $1$ or $\infty$ will be presented in Appendix.  $\hfill\Box$\\

\begin{prop}\label{density} Let $s\in \mathbb{R}$, $1\leq p,q < \infty$.
Then $ \mathscr{S}_m$ is dense in $E^s_{p,q}$.
\end{prop}
{\bf Proof.}  See Appendix. $\hfill\Box$

\subsection{Embeddings between $E^s_{p,q}$ and $B^\sigma_{p,q}$}

Now we give a comparison between $B^s_{p,q}$ and $E^s_{p,q}$, which indicates that $E^{s_1}_{p,q}$ is rougher than all Besov spaces $B^{s_0}_{p,q}$ if $s_1<0$. This property implies that $E^{s}_{p,q}$ with $s<0$  are super-critical spaces for NS.

\begin{prop}\label{compar} Let $s_i\in \mathbb{R}$, $1\leq p_i,q_i \leq  \infty$, $i=0,1$ and $s_1<0$, $p_0\leq p_1$. Then we have
$$
B^{s_0}_{p_0,q_0} \subset  E^{s_1}_{p_1,q_1}.
$$
\end{prop}
{\bf Proof.} By the inclusion $\ell^{q}\subset \ell^r$ for $q\leq r$, it suffices to show that
$$
B^{s_0}_{p_0,\infty} \subset  E^{s_1}_{p_1, 1}.
$$
Let
$$
\Lambda_0=  \{k: \ |k|< 2 \}, \ \   \Lambda_j = \{k: \ |k|\in [2^j, 2^{j+1}) \}, \ \ j\geq 1.
$$
We have
$$
\|f\|_{E^{s_1}_{p_1, 1}} = \sum^\infty_{j=0} \sum_{k\in \Lambda_j} 2^{s_1 |k|}\|\Box_k f\|_{p_1}.
$$
Denote $\tilde{\triangle}_j :=  \triangle_{j-1} + \triangle_{j} + \triangle_{j+1}$. We have for any $k\in \Lambda_j$,
$$
\|\Box_k f\|_p = \|\tilde{\triangle}_j \Box_k f\|_p \leq \|\mathscr{F}^{-1}\sigma_k\|_1  \|\tilde{\triangle}_j  f\|_p \lesssim  \|\tilde{\triangle}_j  f\|_p.
$$
Noticing that $\# \Lambda_j \sim 2^{jd} \sim |k|^d$ for any $k\in \Lambda_j$, we see that
\begin{align}
\|f\|_{E^{s_1}_{p_1, 1}} & \lesssim \sum_{|k|\lesssim 1} \|\Box_k f\|_p +  \sum^\infty_{j=2} \sup_{k\in \Lambda_j} |k|^d 2^{s_1 |k|}\|\Box_k f\|_{p_1} \nonumber\\
& \lesssim \sum_{|k|\lesssim 1} \|\Box_k f\|_p + \sum^\infty_{j=2} 2^{-\varepsilon j} \sup_{k\in \Lambda_j} |k|^{d+\varepsilon } 2^{s_1 |k|}\|\Box_k f\|_{p_0} \nonumber\\
& \lesssim  \sum_{|k|\lesssim 1} \|\Box_k f\|_p +  \| f\|_{B^{s_0}_{p_0, \infty}} \sum^\infty_{j=1} 2^{-\varepsilon j} \sup_{k\in \Lambda_j} |k|^{d+\varepsilon -s_0} 2^{s_1 |k|}\nonumber\\
& \lesssim \| f\|_{B^{s_0}_{p_0, \infty}},
\end{align}
where we have used that $\sup_{k} |k|^\lambda 2^{s_1|k|} <\infty$.  $\hfill\Box$

\begin{prop}\label{embeddings1} Let   $1\leq p, q_i\leq \infty$, $i=1,2$,  $s\in \mathbb{R}$. Then  we have
\begin{align}
& E^{s_+}_{p,q_1}  \subset  E^{s}_{p,q_2},  \label{12embeddings1}\\
& E^{s}_{p,q_1}  \subset  E^{s}_{p,q_2}, \ \ (q_1 \leq q_2),  \label{12embeddings2}\\
& \|2^{s |\xi|} \hat{f}\|_2 \sim \|f\|_{E^s_{2,2}}.
\end{align}
\end{prop}
{\bf Proof.} See \cite{Groch01}. Or, using the embeddings of $\ell^{s}_q$, we can easily get the results, as desired. $\hfill\Box$

\section{Estimates for the linear heat equation}

Recall that the linear heat equation
\begin{align}
u_t -\Delta u = f, \ \ u(0,x) =u_0
\end{align}
has the solution
\begin{align}
u = H(t) u_0 + \mathscr{A} f.
\end{align}
We want to obtain the time-global estimates of $H(t)u_0$ and $\mathscr{A} f$. Observing the symbol $e^{-t|\xi|^2}$ of $H(t)$, which has an exponential decay for $|\xi|> 0$, one can easily get the time-global estimates for the case $|\xi| \gtrsim 1$. But for $|\xi|$ near $0$, the exponential decay of $e^{-t|\xi|^2}$ will heavily rely upon $|\xi|$ and $t>0$, one must carefully treat the dependence of $|\xi|$ to the time-global estimates.

\begin{lem}\label{Groupest1} Let   $1\leq p \leq \infty$, $|k|_\infty\geq 1$. Then  we have
\begin{align}
\|\Box_k H(t)u_0\|_p   \lesssim e^{-ct|k|^2}   \|\Box_{k}  u_0\|_p.  \label{group1}
\end{align}
\end{lem}

{\bf Proof.} \eqref{group1} was essentially obtained  in \cite{Wa06}. We can choose a smooth cut-off function $\rho$ such that ${\rm supp}\, \rho \subset [-7/8,7/8]^d$ and $\rho=1$ in the cube $[-3/4,3/4]^d$. Denote $\rho_k=\rho(\cdot-k)$.  We remark that ${\rm supp}\, \rho_k \subset \{\xi: \ |k_j|-7/8 \leq |\xi_j| \leq \ |k_j|+7/8 \}$ and $\rho_k \sigma_k =\sigma_k$. Hence, $|\xi|\sim |k| \gtrsim 1$ if $|k|_\infty \geq 1$ and $\xi \in {\rm supp}\,\rho_k$. Using Young's and Bernstein's estimates:
\begin{align}
\|\Box_k H(t)u_0\|_p   & \lesssim \|\mathscr{F}^{-1}(\rho_k e^{-t|\xi|^2})\|_1    \|\Box_{k}  u_0\|_p \nonumber\\
 & \lesssim \| \rho e^{-t|\xi+k|^2} \|_{H^{[d/2]+1}}   \|\Box_{k}  u_0\|_p  \nonumber\\
& \lesssim e^{-ct|k|^2}      \|\Box_{k}  u_0\|_p,  \nonumber
\end{align}
which implies the result, as desired. $\hfill\Box$

As a straightforward consequence of Lemma  \ref{Groupest1}, we have

\begin{lem}\label{Groupest3} Let   $1\leq p\leq \infty$, $1\leq \gamma_1\leq \gamma\leq \infty$,  $|k|_\infty\geq 1$. Then  we have
\begin{align}
\|\Box_k H(t)u_0\|_{L^\gamma_t L^p_x}   & \lesssim |k|^{-2/\gamma}   \|\Box_{k}  u_0\|_p.  \label{2group1} \\
\|\Box_k \mathscr{A} f \|_{L^\gamma_t L^p_x}  & \lesssim |k|^{-2(1+ 1 /\gamma-1/\gamma_1)}  \|\Box_{k} f\|_{L^{\gamma_1}_t L^p_x} .  \label{3group1}
\end{align}
\end{lem}

{\bf Proof.} Taking the $L^\gamma_t$ norm in both sides of \eqref{group1} and noticing that $\|e^{-ct|k|^2}\|_{L^\gamma_t} \lesssim  |k|^{-2/\gamma}$, we immediately obtain \eqref{2group1}. By Lemma \ref{Groupest1},
\begin{align}
\|\Box_k \mathscr{A} f \|_{p}   \lesssim \int^t_0  e^{-c(t-\tau)|k|^2}  \|\Box_{k} f(\tau)\|_{p}d\tau .  \label{3group2}
\end{align}
Let $1/\rho = 1+ 1/\gamma-1/\gamma_1$. In view of $\gamma_1\leq \gamma$ we see that $\rho\in [1,\infty]$. By Young's inequality,
\begin{align}
\|\Box_k \mathscr{A} f \|_{L^\gamma_t L^p_x}   \lesssim  \|e^{-c t|k|^2}\|_{L^\rho_t}  \|\Box_{k} f\|_{L^{\gamma_1}_t L^p_x} .  \label{3group3}
\end{align}
Noticing that $\|e^{-t|k|^2}\|_{L^\rho_t} \leq |k|^{-2/\rho}$,  from  \eqref{3group3} we have \eqref{3group1}. $\hfill\Box$

Now we consider the estimate of $\Box_0 H(t) u_0$.  We apply the dyadic decomposition.
In \cite{Che99}, Chemin obtained the following estimate:

\begin{lem}\label{Groupest4} Let   $1\leq  p \leq \infty$. Then  we have
\begin{align}
\|\triangle_j H(t)u_0\|_p   \lesssim e^{-ct 2^{2j}}   \|\triangle_j  u_0\|_p.  \label{4group1}
\end{align}
\end{lem}

\begin{lem}\label{Groupest5} Let   $1 < r\leq  p < \infty$,  $1\leq \gamma \leq \infty$. We have the following results.

\begin{itemize}
\item[\rm (i)] If $ \alpha+ d(1/r-1/p)>2/\gamma$, then we have
\begin{align}
\|(-\Delta)^{\alpha/2} \Box_0 H(t)u_0\|_{L^\gamma_t L^p_x}   \lesssim   \sum_{|\ell| \lesssim 1}  \|\Box_{\ell}  u_0\|_r.  \label{5group2}
\end{align}

\item[\rm (ii)]  Assume that one of the following alternative conditions is satisfied:

\begin{itemize}
\item[\rm (a)] $r\leq \gamma $  and $p>r$;

\item[\rm (b)] $r\vee 2 < \gamma $, $p= r$.
 \end{itemize}

 Then we have
 \begin{align}
\| (-\Delta)^{\alpha/2}  H(t)u_0 \|_{L^\gamma_t L^p_x }    \lesssim     \|   u_0\|_{\dot H^{\alpha  +d(1/r-1/p)-2/\gamma}_{r}}.  \label{5group1a}
\end{align}
In particular, if in addition that $ \alpha+ d(1/r-1/p)=2/\gamma$,   then
\begin{align}
\|(-\Delta)^{\alpha/2}\Box_0 H(t)u_0\|_{L^\gamma_t L^p_x}   \lesssim     \|\Box_{0}  u_0\|_r.  \label{5group1}
\end{align}
 \end{itemize}
 \end{lem}
{\bf Proof.}
By Lemma \ref{Groupest4} and Bernstein's estimates,
\begin{align}
\|(-\Delta)^{\alpha/2} \triangle_j    H(t)  u_0\|_p    \lesssim  e^{-ct 2^{2j}}  2^{\alpha j+ d(1/r-1/p)j} \| \triangle_j     u_0\|_r .  \label{5group4}
\end{align}
Taking $L^\gamma_t$ norm in both sides of \eqref{5group4} and noticing that $\|e^{-ct 2^{2j}}\|_{L^\gamma_t} \lesssim 2^{-2j/\gamma}$,
\begin{align}
\| (-\Delta)^{\alpha/2} \triangle_j    H(t)  u_0\|_{L^\gamma_t L^p_x}   \lesssim    2^{\alpha j -2j/\gamma + d(1/r-1/p)j} \| \triangle_j     u_0\|_r .  \label{5group4a}
\end{align}

 If $\alpha+  d(1/r-1/p)>2/\gamma$, By the definition of $\triangle_j, $ we see that there exists $j_0 \in \mathbb{N}$ such that
$$
\Box_0 = \sum^{j_0}_{j= -\infty} \triangle_j  \Box_0 \ \ {\rm in} \ L^p.
$$
By Lemma \ref{Groupest4} and $\Box_0: L^p \to L^p$,
\begin{align}
\|(-\Delta)^{\alpha/2} \Box_0 H(t)u_0\|_{L^\gamma_t L^p_x}   \leq  \sum^{j_0}_{j= -\infty}  \| \triangle_j  (-\Delta)^{\alpha/2}  H(t)  u_0\|_{L^\gamma_t L^p_x}.  \label{5group3}
\end{align}
it follows from \eqref{5group4a} that
 \begin{align}
\|(-\Delta)^{\alpha/2} \Box_0 H(t)u_0\|_{L^\gamma_t L^p_x}   \lesssim \sup_{j\leq j_0}     \|\triangle_j  u_0\|_r.  \label{5group6}
\end{align}
Noticing that $\sum_{|\ell|_\infty \lesssim 1}  \Box_{\ell} \triangle_j = \triangle_j$ for $j\leq j_0$ and $\triangle_j: L^r \to L^r$, we have
 \begin{align}
      \|\triangle_j  u_0\|_r \lesssim \sum_{|\ell|_\infty \lesssim 1} \|\Box_{\ell} u_0\|_r, \ \ j\leq j_0.  \label{5group7}
\end{align}
In view of \eqref{5group6} and \eqref{5group7} we have \eqref{5group2}.

Now we consider the proof of (ii). If $\gamma =\infty$, noticing that $e^{t\Delta}: L^p\to L^p$ and $\dot H^{d(1/r-1/p)}_r \subset L^p$, we immediately have \eqref{5group1a} for $\gamma=\infty$.  Taking sequence $\ell^q$ norm over all $j$ in both sides of \eqref{5group4a} and using Minkowski's inequality, one has that for any $q \leq \gamma$,
\begin{align}
\| (-\Delta)^{\alpha/2}  H(t)u_0 \|_{L^\gamma_t(\dot B^{2/\gamma}_{p,q})}    \lesssim     \|   u_0\|_{\dot B^{\alpha  +d(1/r-1/p)}_{r,q}}.  \label{5group5}
\end{align}
If $\gamma \geq r, p>r$, then we can choose $p_1, r_1$ satisfying $  p>p_1>r_1>r$. Replacing $r,p$ with $r_1, p_1$ respectively, and taking $q=r$ in \eqref{5group5}, we have
\begin{align}
\| (-\Delta)^{\alpha/2} H(t)u_0 \|_{L^\gamma_t(\dot B^{2/\gamma}_{p_1,r})}    \lesssim     \|   u_0\|_{\dot B^{\alpha  +d(1/r_1-1/p_1)}_{r_1,r}}.  \label{5group8}
\end{align}
 By the embedding $\dot H^{\alpha  +d(1/r  -1/p_1)}_r =\dot F^{\alpha  +d(1/r  -1/p_1)}_{r  ,2} \subset \dot F^{\alpha  +d(1/r  -1/p_1)}_{r  ,\infty} \subset \dot B^{\alpha  +d(1/r_1-1/p_1)}_{r_1,r}$, we have
\begin{align}
\| (-\Delta)^{\alpha/2} H(t)u_0 \|_{L^\gamma_t(\dot B^{2/\gamma}_{p_1,r})}    \lesssim     \|   u_0\|_{\dot H^{\alpha  +d(1/r-1/p_1)}_{r}}.  \label{5group9}
\end{align}
In view of the embedding $\dot B^{2/\gamma}_{p_1,r} \subset \dot B^{2/\gamma}_{p_1,p_1} = \dot F^{2/\gamma}_{p_1,p_1} \subset \dot F^{2/\gamma -d(1/p_1-1/p)}_{p ,2} = \dot H^{2/\gamma  -d(1/p_1-1/p)}_{p}$, we have
\begin{align}
\| (-\Delta)^{\alpha/2} H(t)u_0 \|_{L^\gamma_t(\dot H^{2/\gamma -d(1/p_1-1/p)}_{p})}    \lesssim     \|   u_0\|_{\dot H^{\alpha  +d(1/r-1/p_1)}_{r}}.  \label{5group9abc}
\end{align}
This implies \eqref{5group1a} if $p>r$.

If $\gamma>r\vee 2$ and $p=r$, it follows from the embedding $\dot H^{\alpha  +d(1/r  -1/p)}_r =\dot F^{\alpha  +d(1/r  -1/p)}_{r  ,2} \subset \dot F^{\alpha  +d(1/r  -1/p)}_{r  ,q} \subset \dot B^{\alpha  +d(1/r -1/p)}_{r ,q}$ for $q=r\vee 2$ and \eqref{5group5} that
\begin{align}
\| (-\Delta)^{\alpha/2}  H(t)u_0 \|_{L^\gamma_t(\dot B^{2/\gamma}_{p,q})}    \lesssim     \|   u_0\|_{\dot H^{\alpha  +d(1/r-1/p)}_{r}}.  \label{5group5aaa}
\end{align}
Choosing $\gamma_i > r$ and $1/\gamma = (1-\theta)/\gamma_0 +  \theta /\gamma_1$, $\theta\in (0,1)$, we have
  \begin{align}
\| (-\Delta)^{\alpha/2}  H(t)u_0 \|_{L^{\gamma_i}_t(\dot B^{2/\gamma_i}_{p,\infty})}    \lesssim     \|   u_0\|_{\dot H^{\alpha  +d(1/r-1/p)}_{r}}, \ \ i=0,1.  \label{5group10}
\end{align}
From the Gagliardo-Nirenberg inequality (see \cite{HaMoOzWa11}), we have $\dot B^{2/\gamma_0}_{p,\infty} \cap \dot B^{2/\gamma_1}_{p,\infty} \subset \dot B^{2/\gamma}_{p,1}$ and
\begin{align}
\|f\|_{\dot B^{2/\gamma}_{p,1}} \lesssim  \|f\|^{1-\theta}_{\dot B^{2/\gamma_0}_{p,\infty}}  \|f\|^{\theta}_{\dot B^{2/\gamma_1}_{p,\infty}}.  \label{5group11}
\end{align}
By H\"older's inequality,
\begin{align}
\| f \|_{L^{\gamma}_t(\dot B^{2/\gamma}_{p,1})}  \leq
\| f \|^{1-\theta}_{L^{\gamma_0}_t(\dot B^{2/\gamma_0}_{p,\infty})}   \| f \|^{ \theta}_{L^{\gamma_1}_t(\dot B^{2/\gamma_1}_{p,\infty})}.  \label{5group12}
\end{align}
Taking $f= (-\Delta)^{\alpha/2}  H(t)u_0$ in  \eqref{5group12} and using \eqref{5group10}, we immediately obtain that
\begin{align}
\|  (-\Delta)^{\alpha/2} H(t)u_0 \|_{L^{\gamma}_t(\dot B^{2/\gamma}_{p,1})}   \lesssim     \|   u_0\|_{\dot H^{\alpha  +d(1/r-1/p)}_{r}}.  \label{5group13}
\end{align}
It follows that \eqref{5group1a} and \eqref{5group1} hold if $p=r$,  $ \gamma >r\vee 2$.

$\hfill\Box$

\begin{lem}\label{Groupest6} Let   $1< p_1 \leq  p < \infty$,  $1\leq \gamma_1\leq \gamma \leq \infty$.  We have the following results.

\begin{itemize}

\item[\rm (i)] Assume that
$$
 \frac{d}{p_1} -\frac{d}{p}   +  \frac{2}{\gamma_1}-\frac{2}{\gamma}  >2-\alpha.
$$
Then  we have
\begin{align}
\|\Box_0 (-\Delta)^{\alpha/2} \mathscr{A} f\|_{L^\gamma_t L^p_x}   \lesssim   \sum_{|\ell| \lesssim 1}  \|\Box_{\ell}  f\|_{L^{\gamma_1}_t L^{p_1}_x}.  \label{6group1}
\end{align}

\item[\rm (ii)]  Assume that one of the following  conditions is satisfied:
\begin{itemize}

\item[\rm (a)] $1 < \gamma_1< \gamma <\infty$;

\item[\rm (b)] $ \gamma =\infty$ and $p> p_1 \vee \gamma_1$.

 \end{itemize}

 Then we have
 \begin{align}
  \| (-\Delta)^{\alpha/2} \mathscr{A} f\|_{L^\gamma  (\dot H^{2/\gamma}_{p}) }
    \lesssim    \| f \|_{L^{\gamma_1} (\dot H^{d(1/p_1-1/p)  + 2/\gamma_1   - 2 +\alpha}_{p_1}) }.
    \label{6group1a}
\end{align}
In particular, if
$$
 \frac{d}{p_1} -\frac{d}{p}   +  \frac{2}{\gamma_1}-\frac{2}{\gamma}  = 2-\alpha,
$$
then we have
\begin{align}
  \| (-\Delta)^{\alpha/2} \mathscr{A} f\|_{L^\gamma_t L^p_x }
    \lesssim    \| f \|_{L^{\gamma_1}_t  L^{p_1}_x}.
    \label{6group1b}
\end{align}

 \end{itemize}
\end{lem}
{\bf Proof.} Using the same way as in the proof of Lemma \ref{Groupest5}, we have some $j_0\in \mathbb{N}$,
\begin{align}
\|\Box_0 (-\Delta)^{\alpha/2} \mathscr{A} f\|_{L^\gamma_t L^p_x}   \lesssim   \sum_{j\leq j_0}  \|\triangle_j (-\Delta)^{\alpha/2} \mathscr{A} f\|_{L^{\gamma}_t L^{p}_x}.  \label{6group2}
\end{align}
It follows from Lemma \ref{Groupest4} and Bernstein's inequality that
\begin{align}
  \|\triangle_j  (-\Delta)^{\alpha/2} \mathscr{A} f\|_{p} & \lesssim \int^t_0  e^{-(t-s)2^{2j}} \|\triangle_j (-\Delta)^{\alpha/2}   f(s)\|_{p} ds \nonumber\\
  & \lesssim \int^t_0  e^{-(t-s)2^{2j}} 2^{d(1/p_1-1/p)j +\alpha j} \|\triangle_j f(s)\|_{p_1} ds.
    \label{6group3}
\end{align}
Applying Young's inequality,  from \eqref{6group3}  we have
\begin{align}
  \|\triangle_j  (-\Delta)^{\alpha/2} \mathscr{A} f\|_{L^\gamma_t L^p_x}
    \lesssim  2^{d(1/p_1-1/p)j +(2/\gamma_1 - 2/\gamma)j -(2-\alpha) j} \|\triangle_j f \|_{L^{\gamma_1}_t L^{p_1}_x}.
    \label{6group4}
\end{align}
If $d(1/p_1-1/p)  + 2/\gamma_1 - 2/\gamma >2-\alpha$, we have
\begin{align}
\|\Box_0 (-\Delta)^{\alpha/2} \mathscr{A} f\|_{L^\gamma_t L^p_x}   \lesssim   \sup_{j\leq j_0}  \|\triangle_j  f\|_{L^{\gamma_1}_t L^{p_1}_x}.  \label{6group5}
\end{align}
Noticing that $\sum_{|\ell|_\infty \lesssim 1}  \Box_{\ell} \triangle_j = \triangle_j$ for $j\leq j_0$ and $\triangle_j: L^{p_1} \to L^{p_1}$, we have
 \begin{align}
      \|\triangle_j f\|_{p_1} \lesssim \sum_{|\ell|_\infty \lesssim 1} \|\Box_{\ell} f\|_{p_1}, \ \ j\leq j_0.  \label{6group6}
\end{align}
Hence, \eqref{6group5} and \eqref{6group6} imply that \eqref{6group1} holds.

Next, we prove (ii). First, we consider the case $1<\gamma_1<\gamma<\infty$.  By \eqref{6group3} and $\sup_{t>0} t^{\delta}/e^t \lesssim 1$,
\begin{align}
  \|\triangle_j  (-\Delta)^{\alpha/2} \mathscr{A} f\|_{p}
  & \lesssim \int^t_0   (t-s)^{-\delta}   2^{d(1/p_1-1/p)j +\alpha j-2\delta j} \|\triangle_j f(s)\|_{p_1} ds.
    \label{6group7}
\end{align}
Taking the sequence $\ell^q$ norms in both sides of \eqref{6group7} and using Minkowski's inequality, one obtain that
\begin{align}
  \|  (-\Delta)^{\alpha/2} \mathscr{A} f\|_{\dot B^0_{p,q}}
  & \lesssim \int^t_0   (t-s)^{-\delta} \|  f(s)\|_{\dot B^{d(1/p_1-1/p)  +\alpha - 2\delta }_{p_1,q}}  ds.
    \label{6group8}
\end{align}
Let $\delta= 1+1/\gamma -1/\gamma_1$. Noticing that $1<\gamma_1<\gamma<\infty$, we see that $\delta\in (0,1)$. Applying the Hardy-Littlewood-Sobolev inequality, one has from \eqref{6group8} that
\begin{align}
  \| (-\Delta)^{\alpha/2} \mathscr{A} f\|_{L^\gamma  (\dot B^{2/\gamma}_{p,q}) }
    \lesssim    \| f \|_{L^{\gamma_1} (\dot B^{d(1/p_1-1/p)  + 2/\gamma_1   - 2 +\alpha}_{p_1, q}) }.
    \label{6group9}
\end{align}
\eqref{6group9} implies that
\begin{align}
  \| (-\Delta)^{\alpha/2} \mathscr{A} f\|_{L^\gamma  (\dot B^{2/\gamma}_{p,\infty}) }
    \lesssim    \| f \|_{L^{\gamma_1} (\dot H^{d(1/p_1-1/p)  + 2/\gamma_1   - 2 +\alpha}_{p_1}) }.
    \label{6group10}
\end{align}
By the Gagliardo-Nirenberg inequality and repeating the procedures as in the proof of \eqref{5group13},   we have
\begin{align}
  \| (-\Delta)^{\alpha/2} \mathscr{A} f\|_{L^\gamma  (\dot B^{2/\gamma}_{p,1}) }
    \lesssim    \| f \|_{L^{\gamma_1} (\dot H^{d(1/p_1-1/p)  + 2/\gamma_1   - 2 +\alpha}_{p_1}) },
    \label{6group11}
\end{align}
which implies \eqref{6group1a}.

Now we consider the case $\gamma=\infty$.  By \eqref{6group3} and Young's inequality, we have
\begin{align}
  \|\triangle_j  (-\Delta)^{\alpha/2} \mathscr{A} f\|_{L^\infty_t L^p_x}
    \lesssim  2^{d(1/p_1-1/p)j + 2 j/\gamma_1 -(2-\alpha) j} \|\triangle_j f \|_{L^{\gamma_1}_t L^{p_1}_x}.
    \label{6group12}
\end{align}
Let $q\geq \gamma_1$. Taking sequence $\ell^q$ norm in both sides of \eqref{6group12} and using Minkowski's inequality, one has that
\begin{align}
  \| (-\Delta)^{\alpha/2} \mathscr{A} f\|_{L^\infty  (\dot B^{0}_{p,q}) }
    \lesssim    \| f \|_{L^{\gamma_1} (\dot B^{d(1/p_1-1/p)  + 2/\gamma_1   - 2 +\alpha}_{p_1, q}) }.
    \label{6group13}
\end{align}
Comparing \eqref{5group5} with \eqref{6group13}, we see that $r$ and $p_1$ have the same role when we use the embeddings between Besov and Triebel spaces. One can repeat the procedures in the proof of \eqref{5group1a} to obtain \eqref{6group1a}. Indeed, taking $p>\tilde{p} >\tilde{p}_1>p_1$, we can assume that $\tilde{p} \geq \gamma_1$. By \eqref{6group13},
\begin{align}
  \| (-\Delta)^{\alpha/2} \mathscr{A} f\|_{L^\infty  (\dot B^{0}_{\tilde{p},\tilde{p}}) }
    \lesssim    \| f \|_{L^{\gamma_1} (\dot B^{d(1/\tilde{p}_1-1/\tilde{p})  + 2/\gamma_1   - 2 +\alpha}_{\tilde{p}_1, \tilde{p}}) }.
    \label{6group14}
\end{align}
Noticing the embedding $\dot B^{0}_{\tilde{p},\tilde{p}} \subset \dot F^{d(1/p-1/\tilde{p})}_{p,2}$ and  $ \dot F^{s+ d(1/p_1-1/\tilde{p}_1)}_{p_1,2} \subset  \dot B^{s}_{\tilde{p}_1,\tilde{p}}$, from \eqref{6group14} we get the result, as desired. $\hfill\Box$

 \begin{cor}\label{Groupest8}
 Let   $2\leq r< p <\infty$.
Then  we have
\begin{align}
& \|\Box_k H(t)u_0\|_{L^2_t L^p_x}   \lesssim |k|^{-1}   \|\Box_{k}  u_0\|_{r},  \ \ |k|_\infty\geq 1, \label{8group1} \\
& \|\Box_0 H(t)u_0\|_{L^2_t L^p_x}   \lesssim  \sum_{|\ell|_\infty \lesssim 1} \|\Box_{\ell} (-\Delta)^{-1/2}  u_0\|_{r},  \label{8group2} \\
& \|\Box_m H(t)u_0\|_{L^\infty_t L^r_x}   \lesssim   \|\Box_{m}  u_0\|_{r},  \ \ m \in \mathbb{Z}^d.  \label{8group3}
\end{align}

\end{cor}
{\bf Proof.} By Lemma \ref{Groupest3} and $\|\Box_k f\|_p \lesssim \|\Box_k f\|_{r}$, we have \eqref{8group1}. \eqref{8group2} is a straightforward consequence of Lemma \ref{Groupest5}. Lemma \ref{Groupest1} implies \eqref{8group3}. $\hfill\Box$

\begin{cor}\label{Groupest9} Let   $2\leq p< \infty$, $1\leq r\leq \infty$. Then  we have
\begin{align}
& \|\Box_k \nabla \mathscr{A} f \|_{L^2_t L^p_x}   \lesssim  \sum_{|\ell|_\infty \lesssim 1} \|\Box_{k+\ell} f\|_{L^{1}_t L^{p/2}_x},  \label{9group1} \\
& \|\Box_k \nabla \mathscr{A} f \|_{L^\infty_t L^r_x}   \lesssim  \sum_{|\ell|_\infty \lesssim 1} \|\Box_{k+\ell} f\|_{L^{1}_t L^{r}_x}.  \label{9group2}
\end{align}
\end{cor}
{\bf Proof.} By Lemma \ref{Groupest3} we have \eqref{9group1} and \eqref{9group2} for $|k|_\infty \geq 1$. Taking $\alpha =1$ and $\alpha=0$, respectively in Lemma \ref{Groupest6}, we have \eqref{9group1} and \eqref{9group2} for $k=0$.

 \begin{cor}\label{Groupest10}
Let $d\geq 2$. We have
\begin{align}
& \| \Box_k H(t)u_0\|_{L^{d+2}_{x,t}}   \lesssim \langle k\rangle^{-2/(d+2)}  \|\Box_k u_0\|_{d},  \label{10group1} \\
 & \| \Box_k \nabla \mathscr{A} f \|_{L^{d+2}_{x,t}}   \lesssim    \|\Box_k f\|_{L^{(d+2)/2}_{x,t}}.  \label{10group2}
\end{align}
\begin{align}
& \| \Box_k H(t)u_0\|_{{L}^{\infty}_t L^d_x}   \lesssim   \| \Box_k u_0\|_{d},  \label{10group3} \\
 & \|\Box_k \nabla \mathscr{A} f \|_{ {L}^{\infty}_t L^d_x }   \lesssim    \| \Box_kf\|_{  {L}^{(d+2)/2}_{x,t}}.  \label{10group4}
\end{align}
\end{cor}
{\bf Proof.} Using a similar way as in the proof of Corollaries \ref{Groupest8} and \ref{Groupest9}, we can obtain the result, as desired. $\hfill\Box$

 \begin{cor}\label{Groupest11}
 Let   $2\leq r< d $. There exists $p:=p(r) \gg 1$, such that for any  $p_1<p $ and $r_1<r$, $k\in \mathbb{Z}^d$,
\begin{align}
& \|\Box_k H(t)u_0\|_{L^\infty_t L^r_x}   \lesssim   \|\Box_{k}   u_0\|_{r},   \label{11group3} \\
& \|\Box_k H(t)u_0\|_{L^2_t L^p_x}   \lesssim  \sum_{|\ell|_\infty \lesssim 1} \langle k\rangle^{-1} \|\Box_{k+\ell}   u_0\|_{r},  \label{11group4} \\
& \|\Box_k \nabla \mathscr{A} f \|_{L^2_t L^p_x}   \lesssim  \sum_{|\ell|_\infty \lesssim 1} \|\Box_{ k+\ell} f\|_{L^{1}_t L^{p_1}_x},   \label{11group1} \\
& \|\Box_k  \nabla \mathscr{A} f \|_{L^\infty_t L^r_x}   \lesssim  \sum_{|\ell|_\infty \lesssim 1} \|\Box_{k+ \ell} f\|_{L^{2}_t L^{r_1}_x}.    \label{11group2}
\end{align}

\end{cor}
{\bf Proof.} By Lemma \ref{Groupest3} we have the above estimates for $|k|_\infty \geq 1$. From Lemma \ref{Groupest5} it follows that \eqref{11group3} and  \eqref{11group4} hold for $k=0$. In view of Lemma \ref{Groupest6} we have \eqref{11group1} and  \eqref{11group2} for $k=0$.  $\hfill\Box$

\section{Global well-posedness in $\dot M^{-1}_{r,1}$}

For convenience, we denote
\begin{align}
      \|u\|_{\widetilde{L}^\gamma(I, M^s_{p,q})} =  \left\|\{\langle k\rangle^s \|\Box_k u\|_{L^\gamma_{t\in I}L^p_x} \}_{k\in \mathbb{Z}^d}\right\|_{\ell^q}.  \label{gwns1}
\end{align}
It is known the following nonlinear mapping estimates in $\widetilde{L}^\gamma(I, M^s_{p,q})$ (see \cite{WaHu07}):

 \begin{lem}\label{Nlest1}
 Let   $1\leq p_1, p_2,  p,  \gamma_1, \gamma_2, \gamma\leq \infty$, $s\geq 0$.    Assume further that
$$
\frac{1}{p} =\frac{1}{p_1} +\frac{1}{p_2}, \ \ \frac{1}{\gamma} = \frac{1}{\gamma_1} +\frac{1}{\gamma_2}.
$$
Then  we have
\begin{align}
 \|u_1 u_2\|_{\widetilde{L}^\gamma(I, M^s_{p,1})} \leq  \prod^2_{i=1}  \|u_i\|_{\widetilde{L}^{\gamma_i}(I, M^s_{p_i,1})}.  \label{gwns2}
\end{align}
\end{lem}
{\bf Proof of Theorem \ref{NSthm1}.} We consider the mapping
\begin{align} \label{gwns2a}
\mathscr{T}:  u(t) \to H (t) u_0 +  \mathscr{A} \mathbb{P}\ {\rm div}\,(u\otimes u)
\end{align}
in the space $\widetilde{L}^2(\mathbb{R}_+, M^{0}_{p,1})$. By Corollary \ref{Groupest8}, for $r<p < 2r$,
\begin{align}
\|\Box_k H(t)u_0\|_{\widetilde{L}^2(\mathbb{R}_+, M^{0}_{p,1}) }   & \lesssim \sum_{|k|_\infty \geq 1} |k|^{-1} \|\Box_{k+\ell}  u_0\|_r + \sum_{|\ell|_\infty  \lesssim 1}   \|(-\Delta)^{-1/2}\Box_{\ell}  u_0\|_{r} \nonumber\\
& \lesssim \|u_0\|_{\dot M^{-1}_{r,1}} \label{gwns3}
\end{align}
By Corollary \ref{Groupest9},
\begin{align}
\|\mathscr{A} \mathbb{P}\ {\rm div}\,(u\otimes u)\|_{\widetilde{L}^2(\mathbb{R}_+, M^{0}_{p,1}) }
\lesssim  \sum_{k\in \mathbb{Z}^d}   \|\Box_k (u\otimes u)\|_{L^1_t L^{p/2}_x}.    \label{gwns4}
\end{align}
It follows from \eqref{gwns4} and Lemma \ref{Nlest1} that
\begin{align}
\|\mathscr{A} \mathbb{P}\ {\rm div}\,(u\otimes u)\|_{\widetilde{L}^2(\mathbb{R}_+, M^{0}_{p,1}) }     \lesssim \| u\|^2_{\widetilde{L}^2(\mathbb{R}_+, M^{0}_{p,1}) }.  \label{gwns5}
\end{align}
Hence, in view of \eqref{gwns3} and \eqref{gwns5},
\begin{align}
\|\mathscr{T} u)\|_{\widetilde{L}^2(\mathbb{R}_+, M^{0}_{p,1}) }     \lesssim \|u_0\|_{\dot M^{-1}_{r,1}} +  \| u\|^2_{\widetilde{L}^2(\mathbb{R}_+, M^{0}_{p,1}) }.  \label{gwns6}
\end{align}
By a standard contraction mapping argument, we can show that \eqref{NS1} has a unique solution $u\in \widetilde{L}^2(\mathbb{R}_+, M^{0}_{p,1})$.
By \eqref{9group2} and Lemma \ref{Nlest1}, for $r \geq p/2$
\begin{align}
\|\mathscr{A} \mathbb{P}\ {\rm div}\,(u\otimes u)\|_{{L}^\infty(\mathbb{R}_+, \dot M^{-1}_{r,1}) }
\lesssim  \sum_{k\in \mathbb{Z}^d}   \|\Box_k (u\otimes u)\|_{L^1_t L^{r}_x} \lesssim  \| u\|^2_{\widetilde{L}^2(\mathbb{R}_+, M^{0}_{p,1}) }.    \label{gwns10}
\end{align}
In view of Corollary \ref{Groupest8}, we have from \eqref{gwns2a} and \eqref{gwns10} that
\begin{align}
\|\mathscr{T} u\|_{ {L}^\infty (\mathbb{R}_+, \dot M^{-1}_{r,1}) }  \leq   \sum_{k\in \mathbb{Z}^d}   \|\Box_k (-\Delta)^{-1/2} \mathscr{T} u \|_{L^\infty_t L^{r}_x}   \lesssim \|u_0\|_{\dot M^{-1}_{r,1}} +  \| u\|^2_{\widetilde{L}^2(\mathbb{R}_+, M^{0}_{p,1}) }.  \label{gwns7}
\end{align}
It follows that the solution $u\in {L}^\infty (\mathbb{R}_+, \dot M^{-1}_{r,1})$. $\hfill\Box$

\section{Scaling property of $E^s_{p,q}$}

The scaling properties of ($\alpha$-)modulation spaces has been obtained in \cite{SuTo07,HaWa14}.  For our purpose we consider the dilation property of $E^s_{p,q}$, $s<0$, which is quite different from ($\alpha$-) modulation spaces. Denote $f_\lambda = f(\lambda \, \cdot)$. We have

\begin{prop} \label{scaling}
Let $s<0$, $1\leq p,q\leq \infty$. Then we have
\begin{align}
\|f_\lambda \|_{E^s_{p,q}} \lesssim \lambda^{-d/p} \|f\|_{E^s_{p,q}}   \label{scaling1}
\end{align}
and the constant omitted in the right hand side of \eqref{scaling1} is independent of $\lambda>1$.
\end{prop}

{\bf Proof.} For convenience, we denote
\begin{align}
\Lambda_m = \left\{ k\in \mathbb{Z}^d: m_j- \frac{1}{2} \leq \frac{k_j}{\lambda} < m_j+ \frac{1}{2}, \ j=1,...,d \right\}   \label{scaling2}
\end{align}
We can rewrite $\|f_\lambda\|_{E^s_{p,q}}$ as
\begin{align}
\|f_\lambda \|^q_{E^s_{p,q}} \sim  \sum_{m\in \mathbb{Z}^d} \sum_{k\in \Lambda_m}  2^{s|k|q}  \|\Box_k f_\lambda\|^q_p.  \label{scaling3}
\end{align}
Since $\Box_k f_\lambda = (\mathscr{F}^{-1}\sigma(\lambda\, \cdot-k) \widehat{f})(\lambda \, \cdot)$,  we have
$$
\|\Box_k f_\lambda\|_p \leq \lambda^{-d/p} \|\mathscr{F}^{-1}\sigma(\lambda\, \cdot-k) \widehat{f}\|_p.
$$
For any $k\in \Lambda_m$, we see that $\sigma_n \sigma (\lambda\, \cdot-k) =0$ if $n\neq m+\ell$, $|\ell|_\infty\leq 1$.  It follows that
\begin{align}
\|\Box_k f_\lambda\|_p & \leq \lambda^{-d/p}  \|\mathscr{F}^{-1}\sigma(\lambda\, \cdot-k) \|_1  \sum_{|\ell|_\infty\leq 1} \|\Box_{m+\ell} f\|_p \nonumber\\
 & \lesssim \lambda^{-d/p}    \sum_{|\ell|_\infty\leq 1} \|\Box_{m+\ell} f\|_p.  \label{scaling4}
\end{align}
Then we have for $\lambda\gg 1$,
\begin{align}
I(\lambda) & : =  \sum_{m\in \mathbb{Z}^d,\, |m|\geq 100d} \ \sum_{k\in \Lambda_m}  2^{s|k|q}  \|\Box_k f_\lambda\|^q_p  \nonumber\\
&  \lesssim \lambda^{-dq/p} \sum_{m\in \mathbb{Z}^d, \, |m|\geq 100d } \  \sum_{|\ell|_\infty\leq 1} \|\Box_{m+\ell} f\|^q_p \sum_{k\in \Lambda_m}  2^{s|k|q}   \nonumber\\
&  \lesssim \lambda^{-dq/p} \sum_{m\in \mathbb{Z}^d,\, |m|\geq 50d}   2^{s|m|\lambda q/4} \|\Box_{m} f\|^q_p  \nonumber\\
&  \lesssim \lambda^{-dq/p}  2^{s\lambda q} \|f\|^q_{E^s_{p,q}}.
\label{scaling5}
\end{align}
If $|m| \leq 100d$ we have
$$
\sum_{k\in \Lambda_m} 2^{s|k|q} \lesssim 1.
$$
It follows that
\begin{align}
II(\lambda) & : =  \sum_{m\in \mathbb{Z}^d,\, |m|< 100d} \ \sum_{k\in \Lambda_m}  2^{s|k|q}  \|\Box_k f_\lambda\|^q_p  \nonumber\\
&  \lesssim \lambda^{-dq/p} \sum_{m\in \mathbb{Z}^d, \, |m|<100d }     \|\Box_{m} f\|^q_p  \nonumber\\
&  \lesssim \lambda^{-dq/p}  \|f\|^q_{E^s_{p,q}}.
\label{scaling6}
\end{align}
By the estimates of $I(\lambda)$ and $II(\lambda)$ as in \eqref{scaling5} and \eqref{scaling6}, we have the result, as desired. $\hfill\Box$\\

Roughly speaking, the scaling in $E^s_{p,q}$ $(s<0)$ is the same as those of  $L^p$, which is independent of $q\in [1,\infty]$.  Recall that for the scaling initial data $u_\lambda(0, \cdot)= \lambda u_0(\lambda\, \cdot)$, we have
$$
\|u_\lambda(0, \cdot) \|_{E^s_{p,q}}  \lesssim \lambda^{1-d/p}  \|u_0 \|_{E^s_{p,q}}.
$$
After making the scaling to $u_0$, the initial data $u_\lambda(0, \cdot)$ can be arbitrarily small if we take $\lambda\gg 1$ and $p<d$. In the case $p<\infty$, we have a stronger result:

\begin{prop} \label{2scaling}
Let $s<0$, $1\leq p<\infty$, $1\leq   q\leq \infty$. Then we have for $f\in E^s_{p,q}$,
\begin{align}
  \|f_\lambda \|_{E^s_{p,q}} = o(\lambda^{-d/p}), \ \  \lambda \to \infty.  \label{2scaling1}
\end{align}
\end{prop}

{\bf Proof.} By the estimate of $I(\lambda)$ in \eqref{scaling5},  it suffices to consider the estimate of $II(\lambda)$ in \eqref{scaling6}. Let us observe that $|k| < 200d \lambda$ in $II(\lambda)$. Take a sufficiently small $\varepsilon>0$. We have
\begin{align}
II(\lambda) &  \leq  \sum_{ |k|< 200d\lambda}    2^{s|k|q}  \|\Box_k f_\lambda\|^q_p  \nonumber\\
&  =   \left(\sum_{\varepsilon \lambda < |k|< 200d\lambda } +  \sum_{ |k|\leq  \varepsilon \lambda} \right)  2^{s|k|q}  \|\Box_{k} f_\lambda \|^q_p  \nonumber\\
&  := III(\lambda) + IV (\lambda).
\label{2scaling6}
\end{align}
Take a smooth cut-off function $\rho$, ${\rm supp}\, \rho \subset [-2, 2]^d$ and $\rho =1$ in $[-1,1]^d$. We can assume that $\lambda>1/\varepsilon$. For any $k$, $|k|\leq \varepsilon \lambda$,
\begin{align}
\|\Box_k f_\lambda\|_p & \leq \lambda^{-d/p}  \|\mathscr{F}^{-1}\sigma(\lambda\, \cdot-k) \|_1    \| \mathscr{F}^{-1}\rho( \cdot / 4\varepsilon) \widehat{f}\|_p \nonumber\\
 & \lesssim \lambda^{-d/p}  \| \mathscr{F}^{-1}\rho( \cdot / 4\varepsilon) \widehat{f}\|_p  .  \label{2scaling4}
\end{align}
It follows that
\begin{align}
IV(\lambda)
 & \lesssim \lambda^{-dq/p}  \| \mathscr{F}^{-1}\rho( \cdot / 4\varepsilon) \widehat{f}\|^q_p  .  \label{2scaling3lambda}
\end{align}
For any $k$, $\varepsilon \lambda \leq |k|<200d\lambda$, one can find some $m \in \mathbb{Z}^d$, $|m|<300d$ such that
\begin{align}
\|\Box_k f_\lambda\|_p & \leq \lambda^{-d/p}  \|\mathscr{F}^{-1}\sigma(\lambda\, \cdot-k) \|_1   \sum_{|\ell|_\infty \leq 1} \| \Box_{m+\ell} f\|_p \nonumber\\
 & \lesssim \lambda^{-d/p}     \sum_{|\ell|_\infty \leq 1} \| \Box_{m+\ell} f\|_p.  \label{2scaling4lambda}
\end{align}
It follows from $\sum_{|k|\geq \varepsilon\lambda}  2^{s|k|q} \lesssim 2^{s\varepsilon \lambda q}$ that
\begin{align}
III(\lambda)
 & \lesssim \lambda^{-dq/p} 2^{s\varepsilon \lambda q}   \sup_{|m| \leq 400d} \| \Box_{m} f\|^q_p .  \label{2scaling3lambda1}
\end{align}
Hence, we have for any $\lambda \gg 1/\varepsilon$,
\begin{align}
\|f_\lambda\|_{E^s_{p,q}} \lesssim    \lambda^{-d/p}  \left(  \| \mathscr{F}^{-1}\rho( \cdot / 4\varepsilon) \widehat{f}\|_p +  2^{s\varepsilon \lambda q} \|f\|_{E^s_{p,q}} \right).
\end{align}
Now we show that $\| \mathscr{F}^{-1}\rho( \cdot / 4\varepsilon) \widehat{f}\|_p$ can be arbitrarily small if we take $\varepsilon \ll 1$. We have
$$
(\mathscr{F}^{-1}\rho( \cdot / 4\varepsilon) \widehat{f}) (x) =  c \varepsilon ^d \left(\int_{|y|\leq R}+ \int_{|y|>R} \right)  (\mathscr{F}^{-1}\rho)(  4\varepsilon (x-y)) \Box_0 f(y) dy := A+B.
$$
By Young's inequality,
$$
B  \lesssim \|\mathscr{F}^{-1}\rho\|_1 \|\Box_0 f\|_{L^p(|\cdot| \geq R)}  \lesssim   \|\Box_0 f\|_{L^p(|\cdot| \geq R)},
$$
$$
A  \lesssim \varepsilon^{d-d/p} \|\mathscr{F}^{-1}\rho\|_p \|\Box_0 f\|_{L^1(|\cdot| \leq R)}  \lesssim (\varepsilon R)^{d-d/p}  \|\Box_0 f\|_{L^p(|\cdot| \leq R)}.
$$
By taking $R\gg 1$ and $\varepsilon \ll 1/R$, we see that $\| \mathscr{F}^{-1}\rho( \cdot / 4\varepsilon) \widehat{f}\|_p$  is sufficiently small. $\hfill\Box$

\begin{rem} \rm
In Proposition \ref{scaling}, we see that \eqref{scaling1} is independent of $f\in E^s_{p,q}$ and $\lambda \gg 1$. However, in the proof of  $\| \mathscr{F}^{-1}\rho( \cdot / 4\varepsilon) \widehat{f}\|_p \to 0 $  as $\varepsilon \to 0$, the convergent speed depends on $f$. So, in Proposition \ref{2scaling},  $\lambda^{d/p}\|f_\lambda \|_{E^s_{p,q}} = o(1)$ depends on $f\in E^s_{p,q}$.

\end{rem}

The scaling of $E^s_{p,q}$ becomes worse in the case $\lambda <1$. For our purpose we only consider the scaling in $E^s_{p,1}$ and the scaling in $E^s_{p,q}$  with $q\neq 1$ can be similarly considered.

\begin{prop} \label{3scaling}
Let $s<0$, $1\leq p \leq \infty$. Then we have
\begin{align}
  \|f_\lambda \|_{E^s_{p,1}} \lesssim  \lambda^{-d/p} \|f  \|_{E^{s \lambda}_{p,1}}, \ \  \lambda <1.  \label{3scaling1}
\end{align}
\end{prop}
{\bf Proof.} For convenience, we denote
\begin{align}
\Lambda_k = \left\{ m \in \mathbb{Z}^d: \frac{k_j}{\lambda}- \frac{2}{\lambda} \leq m_j < \frac{k_j}{\lambda} + \frac{2}{\lambda}, \ j=1,...,d \right\}.    \label{3scaling2}
\end{align}
\begin{align}
\|f_\lambda \|_{E^s_{p,1}} \sim  \sum_{k \in \mathbb{Z}^d}    2^{s|k|}  \|\Box_k f_\lambda\|_p.  \label{3scaling3}
\end{align}
Since $\lambda <1$,  we have
\begin{align}
\|\Box_k f_\lambda\|_p & \leq \lambda^{-d/p} \|\mathscr{F}^{-1}\sigma(\lambda\, \cdot-k) \widehat{f}\|_p \nonumber\\
& \lesssim \lambda^{-d/p}  \sum_{m\in \Lambda_k} \|\Box_m f\|_p.
\end{align}
Noticing that $|k| \geq \lambda |m| -C$ if $m\in \Lambda_k$, we have the result, as desired. $\hfill\Box$

\section{Initial data in $E^s_{r,1}, \, s<0$} \label{initialdataneg}

The nonlinear estimates in Besov spaces $B^s_{p,q}$ with $s<0$ seem to be subtle, up to now we could not find any algebraic structure of $B^s_{p,q}$ in the case $s<0$, for instance, one has that
\begin{align}
\|f g \|_{\dot B^{s}_{\infty,\infty}} \lesssim (\|f\|_{L^{\infty}} + \|f\|_{\dot B^{-s}_{-d/s,1}}) \|g\|_{\dot B^{s}_{\infty,\infty}},  \ s<0 \label{nonBesov}
\end{align}
and it is impossible to reduce the regularity of $f$ (below $L^\infty$) to obtain \eqref{nonBesov}, even if $f$ and $g$ satisfy some additional conditions.  Similarly for $M^s_{p,q}$ in the case $s<0$.

However, $E^s_{p,1}$ in the case $s<0$ has some algebraic structure if we consider a class of distributions whose Fourier transforms are ``supported" in one octant.

First, we give some explanation to the support of $f\in E^s_{p,q}$, $s<0$. Noticing that $\Box_k f \in L^p\subset \mathscr{S}'$, one sees that $\chi_{Q_k}\widehat{f} \in \mathscr{S}'$. Recall that for any $f\in \mathscr{S}'$,  ${\rm supp}\, f$ is defined as the smallest closed set $F$ so that $f$ vanishes in the complementary set $F^c$ (i.e., $\langle f, \varphi\rangle=0$ for any test function $\varphi$ with ${\rm supp}\,\varphi \subset F^c$).  It follows that ${\rm supp}\, (\chi_{Q_k}\widehat{f})$ is well-defined. Now we can define
$$
{\rm supp}\, \widehat{f} = \mbox{\rm the   closure    of}  \ \bigcup_{k} {\rm supp}\, (\chi_{Q_k} \widehat{f}), \ \ f\in E^s_{p,q}.
$$
Now we give a nonlinear mapping estimate in $\widetilde{L}^{\gamma} (\mathbb{R}, E^s_{p,1})$. We have

\begin{lem}  \label{7Nlest1}
Let $s \in \mathbb{R}$, $1\leq \gamma, \gamma_1, \gamma_2\leq \infty$,  $1/\gamma = 1/\gamma_1 +1/\gamma_2$ and $1<p,p_1,p_2<\infty$, $1/p = 1/p_1 +1/p_2$. Assume that ${\rm supp}\, \widehat{u}$, ${\rm supp}\, \widehat{v}\subset \mathbb{R}^d_I := \{\xi: \ \xi_i \geq 0, \ i=1,...,d\}$. Then we have
\begin{align}
\| uv \|_{\widetilde{L}^{\gamma}(\mathbb{R}_+ ; E^{s}_{p,1})}  \lesssim  2^{C|s|}  \| u  \|_{\widetilde{L}^{\gamma_1}(\mathbb{R}_+ ; E^{s}_{p_1,1})}  \|v\|_{\widetilde{L}^{\gamma_2}(\mathbb{R}_+ ; E^{s}_{p_2,1})}. \label{step1nonlin}
\end{align}
In particular,
\begin{align} \label{0astep1nonlin}
\| uv \|_{\widetilde{L}^{(d+2)/2}(\mathbb{R}_+ ; E^{s}_{(d+2)/2,1})}  \lesssim 2^{C|s|} \| u  \|_{\widetilde{L}^{d+2}(\mathbb{R}_+ ; E^{s}_{d+2,1})}  \|v\|_{\widetilde{L}^{d+2}(\mathbb{R}_+ ; E^{s}_{d+2,1})}.
\end{align}
\end{lem}
{\bf Proof.}  Let us consider the frequency uniform decomposition of $u,v$, we have for $\Box^c_k = \mathscr{F}^{-1} \chi_{k+[0,1)^d}\mathscr{F}$,
$$
u= \sum_{k\in \mathbb{Z}^d_{I} } \Box^c_k u, \ \  v= \sum_{k\in \mathbb{Z}^d_{I} }  \Box^c_k v, \ \ \mathbb{Z}^d_{I} = \mathbb{R}^d_{I}\cap \mathbb{Z}^d.
$$
Following the proof of Proposition \ref{equivnorm1}, we see that
$$
\| f \|_{\widetilde{L}^{\gamma}(\mathbb{R}_+ ; E^{s}_{p,1})} \sim_s  \sum_{k \in  \mathbb{Z}^d } 2^{s|k|}\|\Box^c_k f \|_{L^\gamma_t L^p_x}
$$
and the equivalent constants in both side can be bounded by $2^{C|s|}$.
Noticing that ${\rm supp}\, \widehat{u}*\widehat{v} \subset \mathbb{R}^d_I$,  we have for  $k, k_1,k_2 \in \mathbb{Z}^d_{I} $,
\begin{align} \label{7nlest1}
2^{s|k|} |\Box^c_k (\Box^c_{k_1} u \Box^c_{k_2} v) |   &= 2^{s|k|} |\Box^c_k (\Box^c_{k_1} u \Box^c_{k_2} v) |  \chi_{\{|k-k_1-k_2|_\infty \leq 3\}}  \nonumber\\
& \leq 2^{3d|s|} 2^{s|k_1|+s|k_2|} | \Box^c_k ( \Box^c_{k_1} u    \Box^c_{k_2} v ) |  \chi_{\{|k-k_1-k_2|_\infty \leq 3\}}
 \end{align}
Taking the summation over all $k, k_1,k_2 \in \mathbb{Z}^d_{I} $ in \eqref{7nlest1}, from Young's inequality we have
\begin{align}
\| uv \|_{\widetilde{L}^{\gamma}(\mathbb{R}_+ ; E^{s}_{p,1})}
& \leq \sum_{k, k_1, k_2 \in \mathbb{Z}^d_{I}  } 2^{C|s|} 2^{s|k|}\|\Box^c_k (\Box^c_{k_1} u  \Box^c_{k_2} v)\|_{L^{p}_{x,t}}  \nonumber\\
& \lesssim  \sum_{k, k_1, k_2 \in \mathbb{Z}^d_{I}  } 2^{C|s|} 2^{s|k_1|+s|k_2|}\|  \Box^c_{k_1} u \|_{L^{\gamma_1}_{t} L^{p_1}_{x}} \| \Box^c_{k_2} v \|_{L^{\gamma_2}_{t} L^{p_2}_{x}} \chi_{\{|k-k_1-k_2|_\infty \leq 3\}} \nonumber\\
& \lesssim 2^{3d|s|} \| u  \|_{\widetilde{L}^{\gamma_1}(\mathbb{R}_+ ; E^{s}_{p_1,1})}  \|v\|_{\widetilde{L}^{\gamma_2}(\mathbb{R}_+ ; E^{s}_{p_2,1})}.
\end{align}
Taking $\gamma=p=(d+2)/2, \ \gamma_1=\gamma_2=\ p_1=p_2=d+2$, we have \eqref{0astep1nonlin}. $\hfill\Box$ \\

One may further ask if the condition ${\rm supp}\, \widehat{u}$, ${\rm supp}\, \widehat{v}\subset \mathbb{R}^d_I$ can be removed in Lemma \ref{7Nlest1}. The answer is definitely negative. In fact, let $u= \Box^c_k (\mathscr{F}^{-1}\sigma_k)$ and $v= \Box^c_{-k} (\mathscr{F}^{-1}\sigma_{-k})$, where $\sigma_k$ is defined in \eqref{sigmak}.  We see that
$$
\|u\|_{E^s_{p,1}} = \|v\|_{E^s_{p,1}} \to 0, \ \ |k|\to \infty.
$$
However, noticing that ${\rm supp}\,\widehat{u}* \widehat{u} \subset [0,2]^d$, we easily see that
$$
\|uv\|_{E^s_{p,1}} \sim 1.
$$

{\bf Proof of Theorem \ref{NSthm2}.} {\it Step} 1. We consider the case $r=d$ and $\|u_0\|_{E^s_{d,1}}$ is sufficiently small. We consider the mapping
\begin{align} \label{initialdataesp11}
\mathscr{T}:  u(t) \to H (t) u_0 +  \mathscr{A} \mathbb{P}\ {\rm div}\,(u\otimes u)
\end{align}
in the space $\widetilde{L}^{d+2}(\mathbb{R}_+, E^{s}_{d+2,1})$ and we introduce
$$
\mathcal{D}_d= \{u\in \widetilde{L}^{d+2}(\mathbb{R}_+, E^{s}_{d+2,1}): \ {\rm supp}\, \widehat{u(t, \cdot)} \subset \mathbb{R}^d_I, \  \|u\|_{\widetilde{L}^{d+2}(\mathbb{R}_+, E^{s}_{d+2,1})} \leq M \}.
$$
For any $f,g$ with ${\rm supp}\,\widehat{f}, \ {\rm supp}\, \widehat{g} \subset \mathbb{R}^d_I$,  we have ${\rm supp}\,(\widehat{f}* \widehat{g}) \subset \mathbb{R}^d_I$. It follows that ${\rm supp}\, \widehat{\mathscr{T}u(t,\cdot)} \subset \mathbb{R}^d_I$ if $\widehat{u(t, \cdot)} \subset \mathbb{R}^d_I$. For any $u\in \mathcal{D} $,  we have
\begin{align} \label{nlmap1}
\|\mathscr{T}  u\|_{\widetilde{L}^{d+2}(\mathbb{R}_+, E^{s}_{d+2,1})} \leq  \|H (t) u_0\|_{\widetilde{L}^{d+2}(\mathbb{R}_+, E^{s}_{d+2,1})} +  \|\mathscr{A} \mathbb{P}\ {\rm div}\,(u\otimes u)\|_{\widetilde{L}^{d+2}(\mathbb{R}_+, E^{s}_{d+2,1})}.
\end{align}
By Corollary \ref{Groupest10},
\begin{align}
   \|H (t) u_0\|_{\widetilde{L}^{d+2}(\mathbb{R}_+, E^{s}_{d+2,1})} \lesssim  \sum_{k} \langle k\rangle^{-2/(2+d)} 2^{s|k|} \|\Box_k u_0\|_d \leq \|u_0\|_{E^s_{d,1}}
\end{align}
and by Corollary \ref{Groupest10} and Lemma \ref{7Nlest1}
\begin{align}\label{nlmap2}
  \|\mathscr{A} \mathbb{P}\ {\rm div}\,(u\otimes u)\|_{\widetilde{L}^{d+2}(\mathbb{R}_+, E^{s}_{d+2,1})} &  \lesssim \| u\otimes u \|_{\widetilde{L}^{(d+2)/2}(\mathbb{R}_+, E^{s}_{(d+2)/2,1})} \nonumber\\
  &  \lesssim \| u  \|^2_{\widetilde{L}^{d+2}(\mathbb{R}_+, E^{s}_{d+2,1})}.
  \end{align}
From \eqref{nlmap1}--\eqref{nlmap2} it follows that
\begin{align} \label{nlmap3}
\|\mathscr{T}  u\|_{\widetilde{L}^{d+2}(\mathbb{R}_+, E^{s}_{d+2,1})}  \lesssim  \|u_0\|_{E^s_{d,1}}  +  \| u  \|^2_{\widetilde{L}^{d+2}(\mathbb{R}_+, E^{s}_{d+2,1})}.
\end{align}
Hence, if $C\|u_0\|_{E^s_{d,1}} \leq M/2$ and $M$ is sufficiently small, we can conclude that $\mathscr{T}: \mathcal{D}\to \mathcal{D}$. Similarly,
$$
\|\mathscr{T}  u - \mathscr{T}  u\|_{\widetilde{L}^{d+2}(\mathbb{R}_+, E^{s}_{d+2,1})} \leq CM  \|  u -   u\|_{\widetilde{L}^{d+2}(\mathbb{R}_+, E^{s}_{d+2,1})},
$$
which implies that $\mathscr{T}$ is contractive in $\mathcal{D}_d$. Again, in view of Corollary \ref{Groupest10}, we see that $u\in C(\mathbb{R}_+, E^s_{d,1})$.

{\it Step } 2. We consider the case $2\leq r<d$ and $\|u_0\|_{E^{s}_{r,1}  }$ is sufficiently small. Let $p\gg 1$  and
$$
\mathcal{D}_r= \{u\in \widetilde{L}^{2}(\mathbb{R}_+, E^{s}_{p,1}): \ {\rm supp}\, \widehat{u(t, \cdot)} \subset \mathbb{R}^d_I, \  \|  u\|_{\widetilde{L}^{2}(\mathbb{R}_+, E^{s}_{p,1})} \leq M \}
$$
and we show that $\mathscr{T}$ is a contraction mapping in $\mathcal{D}_r$. By Corollary   \ref{Groupest11},
\begin{align}
\| H(t)u_0\|_{\widetilde{L}^2(\mathbb{R}_+, E^{s}_{p,1}) }   & \lesssim \|u_0\|_{ E^{s}_{r,1}},  \label{step2gwns3}\\
\|\mathscr{A} \mathbb{P}\ {\rm div}\,(u\otimes u)\|_{\widetilde{L}^2(\mathbb{R}_+, E^{s}_{p,1}) }
& \lesssim  \sum_{k\in \mathbb{Z}^d_I} 2^{|k|s}   \|\Box_k   (u\otimes u)\|_{L^1_t L^{p/2}_x}.    \label{step2gwns4}
\end{align}
It follows from   \eqref{step1nonlin} that
\begin{align}
\|  \mathscr{A} \mathbb{P}\ {\rm div}\,(u\otimes u)\|_{\widetilde{L}^2(\mathbb{R}_+, E^{s}_{p,1}  ) }     \lesssim \|   u\|^2_{\widetilde{L}^2(\mathbb{R}_+, E^{s}_{p,1}) }.  \label{step2gwns5}
\end{align}
Hence, in view of \eqref{step2gwns3} and \eqref{step2gwns5},
\begin{align}
\|  \mathscr{T} u \|_{\widetilde{L}^2(\mathbb{R}_+, E^{s}_{p,1}) }     \lesssim \|u_0\|_{E^{s}_{r,1}  } +  \| u\|^2_{\widetilde{L}^2(\mathbb{R}_+, E^{s}_{p,1} ) }.  \label{step2gwns6}
\end{align}
By a standard contraction mapping argument, we can show that NS has a unique mild solution $u$ with $  u\in \widetilde{L}^2(\mathbb{R}_+, E^{s}_{p,1})$ if $\|u_0\|_{E^{s}_{r,1}  }$ is sufficiently small. Moreover,
\begin{align}
\|    u \|_{\widetilde{L}^2(\mathbb{R}_+, E^{s}_{p,1}) }     \lesssim \|u_0\|_{E^{s}_{r,1}  }  .  \label{1step2gwns6}
\end{align}
By Corollary \ref{Groupest11}, by taking $1/r_1= 1/r+ 1/p$, from \eqref{step1nonlin} we have
\begin{align}
\|\mathscr{A} \mathbb{P}\ {\rm div}\,(u\otimes u)\|_{{L}^\infty(\mathbb{R}_+, E^{s}_{r,1}) }
& \lesssim  \sum_{k\in \mathbb{Z}^d}  2^{|k|s} \|\Box_k (u\otimes u)\|_{L^2_t L^{r_1}_x}  \nonumber\\
& \lesssim  \| u\|_{\widetilde{L}^2(\mathbb{R}_+, E^{s}_{p,1}) } \| u\|_{\widetilde{L}^\infty (\mathbb{R}_+, E^{s}_{r,1}) }.    \label{step2gwns10}
\end{align}
It follows that the solution $u\in C (\mathbb{R}_+, E^{s}_{r,1})$.

{\it Step} 3. We consider any initial value $u_0 \in E^s_{r,1}$. Make a scaling to $u_0$ by letting
$$
u_{\lambda,0}(x) = \lambda u_0(\lambda x), \ \  \lambda \gg 1.
$$
Since $2\leq r\leq d$, in view of Proposition \ref{2scaling} we see that
$$
\|u_{\lambda,0}\|_{E^{s}_{r,1}} \ll 1, \ \ \lambda \gg 1.
$$
By treating $u_{\lambda,0}$ as initial data for NS, we see that NS has a unique mild solution $u_\lambda \in C (\mathbb{R}_+, E^{s}_{r,1}) \cap D_r$. One sees that $u(t,x)= \lambda^{-1} u_\lambda ( \lambda^{-2}t, \lambda^{-1} x)$ is the solution of NS with initial data $u_0$. By Proposition \ref{3scaling}, we see that $u  \in C (\mathbb{R}_+, E^{s \lambda}_{r,1}) \cap \widetilde{L}^{\gamma}(0,\infty; E^{s \lambda}_{p, 1})$ for some $\lambda>1$, where $\gamma=p=2+d$ for $r=d$, $\gamma=2$ and $p\gg 1$ for $2\leq r<d$. $\hfill\Box$

In the proof of Theorem \ref{NSthm2}, the condition for initial data $u_0\in E^s_{r,1}$ can be replaced by the following slightly weaker version:
$$
\sum_{k\in \mathbb{Z}^d} \langle k\rangle^{-2/\gamma} 2^{s|k|} \|\Box_k u_0\|_r  <\infty,
$$
where $\gamma=2+d$ for $r=d$, $\gamma=2$ for $2\leq r<d$. Using the embedding $E^{s_+}_{2,2} \subset E^s_{2,1} \subset E^s_{r,1}$ for any $r>2$, we see that $u_0\in E^s_{r,1}$  can be replaced by $u_0\in E^{s_+}_{2,2}$, $\|u_0\|_{E^{s_+}_{2,2}} \sim \|2^{s|\cdot|} \widehat{f} \|_2$.

\section{Analyticity of solutions}

In this subsection, we show that for the sufficiently small initial data $u_0\in E^{s}_{r,1}$ with ${\rm supp}\,\widehat{u}_0 \subset \mathbb{R}^d_I$, the corresponding solution will be analytic after the time $t = (s/c)^2$. Since the scaling solution $u_\lambda (t,x) = \lambda u(\lambda^2 t, \lambda x)$ with initial data $u_{0,\lambda}$ which are sufficiently small in $E^s_{r,1}$, we obtain that there exists $\lambda_0 \gg 1$ such that the solution $u_{\lambda_0}$ has the analyticity after the time $t = (s/c)^2$. It follows that $u$ has the analyticity after $t=(\lambda_0 s)^2/c^2$.

Let $|k|_\infty\geq 1$. By Lemma \ref{Groupest1}, we see that
\begin{align}
2^{c\sqrt{t} |k|/2}\|\Box_k H(t)u_0\|_p   \lesssim e^{-ct|k|^2/2}   \|\Box_{k}  u_0\|_p.  \label{smgroup1}
\end{align}
Taking the $L^\gamma_t$ norms in both sides of \eqref{smgroup1}, we immediately have

\begin{lem}\label{smGroupest2}
Let   $1\leq p, \gamma \leq \infty$, $|k|_\infty\geq 1$. Then there  exists $c>0$ such that
\begin{align}
\|2^{c\sqrt{t} |k|}\Box_k H(t)u_0\|_{L^\gamma_t L^p_x}   \lesssim |k|^{-2/\gamma}  \|\Box_{k}  u_0\|_p.  \label{smgroup2}
\end{align}
\end{lem}
Again, in view of Lemma \ref{Groupest1},
\begin{align}
\|2^{c\sqrt{t} |k|/2} \Box_k \mathscr{A} f \|_{p}   \lesssim \int^t_0  e^{-c(t-\tau)|k|^2/2}   \|2^{c\sqrt{\tau} |k|/2}\Box_{k} f(\tau)\|_{p}d\tau .  \label{smgroup3}
\end{align}
By Young's inequality, we have

\begin{lem}\label{smGroupest3} Let   $1\leq p\leq \infty$, $1\leq \gamma_1\leq \gamma\leq \infty$,  $|k|_\infty\geq 1$. Then  there exists $c>0$ such that
\begin{align}
\|2^{c\sqrt{t} |k|} \Box_k \mathscr{A} f \|_{L^\gamma_t L^p_x}   \lesssim |k|^{-2(1+ 1 /\gamma-1/\gamma_1)} \|2^{c\sqrt{t} |k|} \Box_{k} f\|_{L^{\gamma_1}_t L^p_x} .  \label{3smgroup1}
\end{align}
\end{lem}

 \begin{cor}\label{smGroupest10}
Let $d\geq 2$. We have for $k\in \mathbb{Z}^d$
\begin{align}
& \|2^{c\sqrt{t} |k|} \Box_k H(t)u_0\|_{L^{d+2}_{x,t}}   \lesssim \langle k\rangle^{-2/(d+2)}  \|\Box_k u_0\|_{d},  \label{sm10group1} \\
 & \|2^{c\sqrt{t} |k|} \Box_k \nabla \mathscr{A} f \|_{L^{d+2}_{x,t}}   \lesssim    \|2^{c\sqrt{t} |k|}\Box_k f\|_{L^{(d+2)/2}_{x,t}}.  \label{sm10group2}
\end{align}
\begin{align}
& \| 2^{c\sqrt{t} |k|}\Box_k H(t)u_0\|_{{L}^{\infty}_t L^d_x}   \lesssim   \| \Box_k u_0\|_{d},  \label{sm10group3} \\
 & \|2^{c\sqrt{t} |k|} \Box_k \nabla \mathscr{A} f \|_{ {L}^{\infty}_t L^d_x }   \lesssim    \|2^{c\sqrt{t} |k|} \Box_kf\|_{  {L}^{(d+2)/2}_{x,t}}.  \label{sm10group4}
\end{align}
\end{cor}
{\bf Proof.} If $|k|_\infty \geq 1$, the result follows from Lemmas \ref{smGroupest2} and \ref{smGroupest3}. If $k=0$, the results have been shown in Corollary \ref{Groupest10}. $\hfill\Box$

Recall that for $s(t) : \mathbb{R}_+\to \mathbb{R}$, we denote
\begin{align}
      \|u\|_{\widetilde{L}^\gamma(I, E^{s(t)}_{p,1})} =  \left\|\{ \|2^{s(t) |k|} \Box_k u\|_{L^\gamma_{t\in I}L^p_x} \}_{k\in \mathbb{Z}^d}\right\|_{\ell^1}.  \label{smNSspace1}
\end{align}
Take $s(t)=s+c\sqrt{t}$. It follows from Corollary \ref{smGroupest10} that
\begin{align}\label{smnlmap1}
  & \|H (t) u_0\|_{\widetilde{L}^{d+2}(\mathbb{R}_+, E^{s+c\sqrt{t}}_{d+2,1})} \lesssim  \sum_{k} \langle k\rangle^{-2/(2+d)} 2^{s|k|} \|\Box_k u_0\|_d \leq \|u_0\|_{E^s_{d,1}} \\
&\|\mathscr{A} \mathbb{P}\ {\rm div}\,(u\otimes u)\|_{\widetilde{L}^{d+2}(0,T; E^{s+c\sqrt{t}}_{d+2,1})}    \lesssim \| u\otimes u \|_{\widetilde{L}^{(d+2)/2}(0,T; E^{s+c\sqrt{t}}_{(d+2)/2,1})} . \label{smnlmap2}
  \end{align}
So, one needs to make a bilinear estimate in $\widetilde{L}^{(d+2)/2}(\mathbb{R}_+, E^{s+c\sqrt{t}}_{(d+2)/2,1})$. We can imitate the ideas in Lemma \ref{Nlest1} to show that

\begin{lem}  \label{smNlest1}
Let $s<0$, ${\rm supp}\, \widehat{u}$, ${\rm supp}\, \widehat{v}\subset \mathbb{R}^d_I := \{\xi: \ \xi_i \geq 0, \ i=1,...,d\}$. Then we have
\begin{align}
\| uv \|_{\widetilde{L}^{(d+2)/2}(0,T ; E^{s+c\sqrt{t}}_{(d+2)/2,1})}  \lesssim 2^{C(|s|+c\sqrt{T})}  \| u  \|_{\widetilde{L}^{d+2}(0,T; E^{s+c\sqrt{t}}_{d+2,1})}  \|v\|_{\widetilde{L}^{d+2}(0,T; E^{s+c\sqrt{t}}_{d+2,1})}
\end{align}
\end{lem}
{\bf Proof.} We use the same notations as in Lemma \ref{Nlest1}.
Following the proof of Proposition \ref{equivnorm1}, we see that
$$
\| f \|_{\widetilde{L}^{\gamma}(0,T ; E^{(s+c\sqrt{t})}_{p,1})} \sim  \sum_{k \in  \mathbb{Z}^d } 2^{(s+c\sqrt{t})|k|}\|\Box^c_k f \|_{L^\gamma (0,T; L^p)},
$$
where the equivalent constants are $2^{\pm C(|s|+c\sqrt{T})}$. Noticing that ${\rm supp}\, \widehat{u}*\widehat{v} \subset \mathbb{R}^d_I$,  we have for  $k, k_1,k_2 \in \mathbb{Z}^d_{I} $,
\begin{align}
2^{(s+c\sqrt{t})|k|} |\Box^c_k (\Box^c_{k_1} u \Box^c_{k_2} v) |   &= 2^{(s+c\sqrt{t})|k|} |\Box^c_k (\Box^c_{k_1} u \Box^c_{k_2} v) |  \chi_{\{|k-k_1-k_2|_\infty \leq 3\}}  \nonumber\\
& \leq 2^{3(|s|+c\sqrt{t})} 2^{(s+c\sqrt{t})|k_1|+(s+c\sqrt{t})|k_2|} | \Box^c_k ( \Box^c_{k_1} u    \Box^c_{k_2} v ) |  \chi_{\{|k-k_1-k_2|_\infty \leq 3\}}. \nonumber
 \end{align}
After showing the dependence of the exact constants to $T$, we can repeat the procedure as in the proof of Lemma \ref{Nlest1} to get the result, as desired. $\hfill\Box$

Now let us fix $T>0$ such that $s+c\sqrt{T}=3|s|$. We can use the same way to construct contraction mapping as in Section \ref{initialdataneg}. Put
$$
\mathcal{D}= \{\widetilde{L}^{d+2}(0,T; E^{s+c\sqrt{t}}_{d+2,1}): \ {\rm supp}\, \widehat{u(t, \cdot)} \subset \mathbb{R}^d_I, \  \|u\|_{\widetilde{L}^{d+2}(0,T; E^{s+c\sqrt{t}}_{d+2,1})} \leq M \}.
$$
We consider the mapping
\begin{align} \label{initialdataesp11a}
\mathscr{T}:  u(t) \to H (t) u_0 +  \mathscr{A} \mathbb{P}\ {\rm div}\,(u\otimes u)
\end{align}
in the space $\mathcal{D}$. We have from \eqref{smnlmap1}, \eqref{smnlmap2} and Lemma \ref{smNlest1} that
\begin{align} \label{1dataespnl}
\|\mathscr{T} u\|_{\widetilde{L}^{d+2}(0,T; E^{s+c\sqrt{t}}_{d+2,1})}  \lesssim  \|u_0\|_{E^s_{d,1}}+ \|  u\|^2_{\widetilde{L}^{d+2}(0,T; E^{s+c\sqrt{t}}_{d+2,1})}.
\end{align}
Then, we can show by a standard contraction mapping argument that NS has a unique solution $u\in \widetilde{L}^{d+2}(0,T; E^{s+c\sqrt{t}}_{d+2,1})$ if $\|u_0\|_{E^s_{d,1}}$ is sufficiently small. Moreover, in view of Corollary \ref{smGroupest10} we see that  $u\in \widetilde{L}^{\infty}(0,T; E^{s+c\sqrt{t}}_{d,1})$.

Recall that $\|u_0\|_{E^s_{d,1}}$ is sufficiently small implies that NS is globally well-posed in $E^s_{d,1}$, the solution obtained in $\widetilde{L}^{d+2}(0,T; E^{s+c\sqrt{t}}_{d+2,1})$ must coincide with the global solution. Now let $s+c\sqrt{t}_0=|s|$, one can extend the solution starting from $t=t_0$ and consider the mapping
\begin{align} \label{initialdataesp12}
\mathscr{T}:  u(t) \to H (t-t_0) u(t_0) +  \int^t_{t_0} e^{(t-\tau)\Delta} \mathbb{P}\ {\rm div}\,(u\otimes u)(\tau) d\tau.
\end{align}
Using the same way as in \eqref{1dataespnl}, we have
\begin{align} \label{1dataespn2}
\|\mathscr{T} u\|_{\widetilde{L}^{d+2}(t_0, t_0+T_1; E^{s+c\sqrt{t-t_0}}_{d+2,1})}  \lesssim  \|u(t_0)\|_{E^s_{d,1}}+ \|  u\|^2_{\widetilde{L}^{d+2}(t_0, t_0+T_1; E^{s+c\sqrt{t-t_0}}_{d+2,1})}.
\end{align}
Noticing that $\|u(t_0)\|_{E^s_{d,1}} \lesssim \|u_0\|_{E^s_{d,1}}$, we can assume that $T_1=T$ and obtain a unique solution $u\in \widetilde{L}^{d+2}(t_0, t_0+T; E^{s+c\sqrt{t-t_0}}_{d+2,1}) \cap \widetilde{L}^{\infty}(t_0, t_0+T; E^{s+c\sqrt{t-t_0}}_{d,1})$.  Repeating the procedures as above, we obtain that the solution of NS $u\in \widetilde{L}^{d+2}( 0, \infty; E^{s(t)}_{d+2,1}) \cap \widetilde{L}^{\infty}( 0, \infty; E^{s(t)}_{d,1})$, $s(t)=\min (s+c\sqrt{t}, |s|)$.

In the case $2\leq r<d$, the argument is similar to the case $r=d$ by using the following result:

\begin{cor}\label{smGroupestr<d}
Let $2\leq r<d$.  There exists  $  p:=p(r) \gg 1$ and $c>0$ such that for any $p_1<p$ and $r_1<r$,
\begin{align}
& \|2^{c\sqrt{t} |k|}\Box_k H(t)u_0\|_{L^2_t L^p_x}   \lesssim \langle k\rangle^{-1}  \|\Box_{k}  u_0\|_r,  \label{smgroup2cor}\\
 & \|2^{c\sqrt{t} |k|} \Box_k \nabla \mathscr{A} f \|_{L^2_t L^p_x}   \lesssim  \|2^{c\sqrt{t} |k|} \Box_{k} f\|_{L^{1}_t L^{p_1}_x} ,  \label{3group1cor}\\
& \|2^{c\sqrt{t} |k|}\Box_k H(t)u_0\|_{L^\infty_t L^r_x}   \lesssim   \|\Box_{k}  u_0\|_r,  \label{smgroup2cor1}\\
 & \|2^{c\sqrt{t} |k|} \Box_k \nabla \mathscr{A} f \|_{L^\infty_t L^r_x}   \lesssim  \|2^{c\sqrt{t} |k|} \Box_{k} f\|_{L^{2}_t L^{r_1}_x} .  \label{3group1cor1}
\end{align}
\end{cor}
{\bf Proof.} The results in the case $k=0$  are the same ones as Corollary \ref{Groupest11} and the case $k\neq 0$ is from Lemmas \ref{smGroupest2} and \ref{smGroupest3}. $\hfill\Box$ \\

\noindent Using the same way as the case $r=d$, we can show that the solution $u$ satisfying
$$
 u\in \widetilde{L}^{2}( 0, \infty; E^{s(t)}_{p,1}) \cap \widetilde{L}^{\infty}( 0, \infty; E^{s(t)}_{r,1}), \ \  s(t)=\min (s+c\sqrt{t}, |s|)
$$
if $u_0 \in E^s_{r,1}$ is sufficiently small and the details are omitted.

If $s(t)=s+c \sqrt{t}>0$,  we see that $E^{s(t)}_{r, 1} \subset G_{1,r}$  (see Appendix), where the functions are analytic. So, the solution $u_{\lambda_0}$ is spatial analytic if $t>s^2/c^2$. It follows that $u$  is spatial analytic if $t>(\lambda_0 s)^2/c^2$.

\section{Appendix}

\subsection{Proofs of some results in Sections 3 and 5}

Following Theorem 11.4.2 in \cite{Groch01}, Propositions \ref{GelfandShilov}, \ref{Banachspace} and \ref{density} can be shown by imitating the arguments as in those of modulation spaces $M^s_{p,q}$  in \cite{Groch01}. Here we give the details of the proofs.

{\bf Proof of Proposition \ref{GelfandShilov}.} For any $f\in \mathscr{S}_m$, we have (cf. \cite{Groch01})
$$
f(x) = c  \int_{\mathbb{R}^{2d}} V_g f (y,\xi) (M_\xi T_y g(x))  dyd\xi.
$$
It follows from $e^{\lambda|x|} \leq e^{\lambda|y|} e^{\lambda|x-y|}$ that
$$
e^{\lambda|x|} |f(x)|  \lesssim  \| f \|_{M^{m_\lambda}_{1,1}}  \sup_{x\in \mathbb{R}^d} e^{\lambda|x|} |g(x)|  \lesssim  \| f \|_{M^{m_\lambda}_{1,1}}.
$$
Hence $p_{\lambda,1}(f) \lesssim  \| f \|_{M^{m_\lambda}_{1,1}} $.  Replacing $f$ by $\widehat{f}$ in the above inequality and noticing that  $\| f \|_{M^{m_\lambda}_{1,1}} = \| \widehat{f} \|_{M^{m_\lambda}_{1,1}}$, we have $q_{\lambda,1}(f) \lesssim  \| f \|_{M^{m_\lambda}_{1,1}} $.

Conversely, for any $f\in \mathscr{S}_1$, we have from \eqref{STFTprop} that
$$
e^{\lambda|\xi|} |V_g f(x,\xi)|  \lesssim  e^{\lambda|\xi|}   \int_{\mathbb{R}^{d}}  e^{-|\xi-\eta|^2/2} |\widehat{f}(\eta)| d\eta  \lesssim q_{\lambda,1} (f) \int_{\mathbb{R}^{d}}    e^{\lambda|\xi|} g(\xi)d\xi   \lesssim q_{\lambda,1}(f)
$$
and
$$
e^{\lambda|x|} |V_g f(x,\xi)|  \lesssim e^{\lambda|x|} |(\mathscr{F}^{-1} g(\xi-\cdot))| * |f(y)|   \lesssim  p_{\lambda,1} (f).
$$
Noticing that
$$
\|f\|_{M^{m_\lambda}_{1,1}} \leq \sup_{x,\xi} (e^{2\lambda|x|} +  e^{2\lambda|\xi|}) |V_g f(x,\xi)|,
$$
we immediately have the result, as desired. $\hfill\Box$\\

\noindent {\bf Proof of Proposition \ref{Banachspace}.} Since $\|\cdot\|_{E^s_{p,q}}$ and  $\|\cdot\|^{\circ}_{E^s_{p,q}}$ are equivalent norms, a Cauchy sequence $\{f_n\} $ in $E^s_{p,q}$ satisfies
$$
\|f_n-f_m\|^{\circ}_{E^s_{p,q}} \to 0 , \ \ m,n \to \infty.
$$
It follows that $\{\Box_kf_n \}$ is a Cauchy sequence in $L^p$, which has a limit $f_k\in L^p$. We may assume that ${\rm supp}\,\widehat{f}_k \subset k+[-3/2,3/2]^d$. In fact, in view of Bernstein's estimate,
$$
\left\|\Box_k f_n-  \sum_{|\ell|_\infty \leq 1} \Box_{k+\ell}f_k \right\|_p \lesssim \|\Box_k f_n-f_k\|_p  \to 0 , \ \  n \to \infty,
$$
one can replace $f_k $ by $\sum_{|\ell|_\infty \leq 1} \Box_{k+\ell}f_k$ if the support of $\widehat{f}_k$ not contained in $k+[-3/2,3/2]^d$.  Put $f=\sum_k f_k$, we see that
$$
f_n -f = \sum_k (\Box_k f_n -f_k).
$$
By Fatou's Lemma, we have
$$
\|f_n-f\|_{E^s_{p,q}}  \leq \|f_n-f_m\|^{\circ}_{E^s_{p,q}} \to 0.
$$
So, $(E^s_{p,q}, \, \|\cdot\|_{E^s_{p,q}})$ is complete.   $\hfill\Box$\\

In order to show the density of $\mathscr{S}_m$ in $E^s_{p,q}$, we need a Gabor expansion in modulation spaces (cf. \cite{Groch01}).

\begin{prop}\label{Gaborexp} Let $s\in \mathbb{R}$, $1\leq p,q< \infty$.
Then $ f\in {M^s_{p,q}}$  has an expansion
$$
f= \sum_{m,l \in \mathbb{Z}^d} c_{ml} e^{imx} \exp(-|x-l|^2/2)
$$
and
$$
\|f\|_{M^s_{p,q}} \sim \left\|\|\{\langle m\rangle^s c_{ml}\}\|_{\ell^p_{l\in \mathbb{Z}^d}} \right\|_{\ell^q_{m\in \mathbb{Z}^d} }.
$$
The same expansion holds for $f\in M^{m_\lambda}_{1,1}$.
\end{prop}

\noindent {\bf Proof of Proposition \ref{density}. } It suffices to consider the case $s<0$.  By the equivalent norm on $E^s_{p,q}$, we see that
$$
\left\|\sum_{|k|\leq K}  \Box_k f -f \right\|^{\circ}_{E^s_{p,q}} \to 0, \ \ K \to \infty.
$$
For any $K\in \mathbb{N}$, we have $F_K:= \sum_{|k|\leq K}  \Box_k f \in M_{p,q}$, which has a Gabor expansion
\begin{align} \label{Gabor}
F_K = \sum_{m,l \in \mathbb{Z}^d} c_{ml} e^{imx} \exp(-|x-l|^2/2) \ \ {\rm in} \ M_{p,q}.
\end{align}
Since $M_{p,q} \subset E^s_{p,q}$, we see that  \eqref{Gabor} also holds in $E^s_{p,q}$.  Noticing that $M_\xi T_x e^{-|\cdot|^2/2} \in \mathscr{S}_m$, we have the result, as desired. $\hfill\Box$\\

\begin{prop}
$\mathscr{S}_1$ is dense in $M^{m_\lambda}_{1,1}$.
\end{prop}
{\bf Proof.} It is a direct consequence of the Gabor expansion in $M^{m_\lambda}_{1,1}$,  cf. \cite{Groch01}, Theorem 13.6.1.
 $\hfill \Box$

\subsection{Relations between $E^s_{p,q}$ and Gevrey class}

Let $\alpha=(\alpha_1, ..., \alpha_d)$,  $\alpha! = \alpha_1!...\alpha_d!$  and $\partial^\alpha=\partial^{\alpha_1}_{x_1}... \partial^{\alpha_n}_{x_d}$.
Recall that the Gevrey class is defined as follows.
$$  G_{1, p}= \left\{f\in C^\infty(\mathbb{R}^d):\,  ^\exists\rho,  M>0\; s.t.\;
 \|\partial^\alpha f \|_p \le   \frac{M \alpha!}{\rho^{|\alpha|}}  ,  \ \forall \
 \alpha\in {\Bbb
Z}^d_+  \right\}.
$$
It is known that $G_{1,\infty}$ is the Gevrey 1-class. Moreover, we easily see that $G_{1,p_1} \subset G_{1, p_2}$ for any $p_1\le p_2$.  There is a very beautiful relation between Gevrey class and exponential modulation spaces, we can show that  $G_{1,p}=\bigcup_{s>0}E^s_{p, q}$.  Roughly speaking, $E^s_{p, q}$ can be regarded as modulation spaces with analytic regularity. It is easy to see that $M^0_{p,q}= E^{0}_{p,q}$.

\begin{prop}\label{prop4.1}
{\rm ({ \cite{Tr83},  Nikol'skij's inequality})} \it Let $\Omega\subset
\mathbb{R}^d$ be a compact set,  $0<r\le \infty.$ Let us denote
$\sigma_r=d( 1/(r\wedge 1)-1/2)$ and assume that $s>\sigma_r$. Then
there exists a constant $C>0$ such that
\begin{align*}
\|\mathscr{F}^{-1}\varphi \mathscr{F} f\|_r\le
C\|\varphi\|_{H^s}\|f\|_r
\end{align*}
holds for all $f\in L^r_\Omega:=\{f\in L^r: \; {\rm supp}\,
\widehat{f} \subset \Omega\}$ and $\varphi\in H^s$. Moreover,  if
$r\ge 1$,  then the above inequality holds for all $f\in L^r$.
\end{prop}

\begin{prop}\label{prop4.2}
Let $ 0<  p ,  q \leq \infty$. Then
\begin{align}
 G_{1, p} = \bigcup_{s>0} {E}^{s}_{p, q}.
\end{align}
\end{prop}

{\it Proof.} We have
\begin{align}
 \|\partial^\alpha f \|_p \le \sum_{k\in \mathbb{Z}^d} \|\Box_k \partial^\alpha f \|_p.
\end{align}
We easily see that
\begin{align}
  \|\Box_k \partial^\alpha f \|_p \le \sum_{|l|_\infty\le 1} \|\sigma_{k+l} \xi^\alpha \|_{M_p} \|\Box_k f\|_p.
\end{align}
Since $M_p$ is translation invariant,  in view of Nikol'skij's inequality we have
\begin{align}
 \|\sigma_{k} \xi^\alpha \|_{M_p} = \|\sigma \cdot (\xi+k)^\alpha \|_{M_p} \lesssim  \|\sigma \cdot (\xi+k)^\alpha \|_{H^L},  \ \ L>(d/(1\wedge p) -1/2).
\end{align}
It is easy to calculate that
\begin{align}
  \|(\xi+k)^\alpha \sigma   \|_{H^L} & \lesssim \|(\xi+k)^\alpha \sigma\|_2 +  \sum^d_{i=1} \|\partial^L_{x_i}((\xi +k)^\alpha \sigma)\|_2 \nonumber\\
& \lesssim \prod^d_{i=1} \langle k_i\rangle^{\alpha_i} + \sum^d_{i=1} \sum^L_{l_i=1}  \prod^{l_i-1}_{m=0} (\alpha_i -m)   \prod^d_{j=1,  j\neq i} \langle k_j\rangle^{ \alpha_j} \langle k_i\rangle^{ \alpha_i-l_i}.
\end{align}
It follows that
\begin{align}
 \|\partial^\alpha f \|_p
 & \lesssim     \sum_{k\in \mathbb{Z}^d}   \prod^d_{i=1} \langle k_i\rangle^{\alpha_i}   \|\Box_k   f \|_p \nonumber\\
  & + \sum^d_{i=1} \sum^L_{l_i=1} \ \sum_{k\in \mathbb{Z}^d} \prod^{l_i-1}_{m=0} (\alpha_i -m)   \prod^d_{j=1,  j\neq i} \langle k_j\rangle^{ \alpha_j} \langle k_i\rangle^{ \alpha_i-l_i}  \|\Box_k   f \|_p. \nonumber
\end{align}
If we can show that
\begin{align} \label{3.8}
  \frac{s^{|\alpha|}} {\alpha !} \frac{\prod^{l_i-1}_{m=0} (\alpha_i -m)   \prod^d_{j=1,  j\neq i} \langle k_j\rangle^{ \alpha_j} \langle k_i\rangle^{ \alpha_i-l_i} }{ e^{2s(|k_1|+...+|k_d|)}} \lesssim 1,
\end{align}
then together with
 $ E^{s+ \varepsilon}_{p, q} \subset E^{s}_{p, 1}$ for all $ s,  \varepsilon>0$ (cf. \cite{Wa06}),   we immediately have
 \begin{align} \label{3.9}
 \|\partial^\alpha f \|_p
 \lesssim   \frac{\alpha !} {s^{|\alpha|}}  \sum_{k\in \mathbb{Z}^d}  e^{2s(|k_1|+...+|k_d|)}     \|\Box_k   f(x)\|_p \lesssim   \frac{\alpha !} {s^{|\alpha|}}       \|f\|_{E^{cs}_{p, 1}} \lesssim   \frac{\alpha !} {s^{|\alpha|}}       \|f\|_{E^{\tilde{c}s}_{p, q}}
\end{align}
for any $ \tilde{c}>c= 2d \mathrm{log}_2 e$.   Noticing that $s>0$ is arbitrary,   it follows from \eqref{3.9} that
\begin{align}
\bigcup_{s>0} {E}^{s}_{p, q} \subset  G_{1, p} .
\end{align}

Now we show \eqref{3.8}.  Using Taylor's expansion,  we see that
\begin{align}
 e^{2s |k_i| } \ge \frac{(2s |k_i|)^{ \alpha_i -l_i  }}{(\alpha_i -l_i) !},  \quad e^{2s |k_j| } \ge \frac{(2s |k_j|)^{ \alpha_j }}{ \alpha_j !},
 \end{align}
 from which we see that \eqref{3.8} holds.

Next,  we show that   $ G_{1, p}  \subset \bigcup_{s>0} E^s_{p,  \infty}$.  We have
\begin{align} \label{3.12}
 \|f\|_{E^{s\mathrm{log}_2 e}_{p,  \infty}} & = \sup_k e^{s|k|} \|\Box_k f\|_p \nonumber\\
 & \lesssim   \sup_k  \sum^\infty _{m=0}  \frac{  (ds)^{ m} }{ m !} (|k_1 |^m+...+ |k_d|^m) \|\Box_k f\|_p
\end{align}
If $|k_i| \le 10$,  we see that
\begin{align} \label{3.13}
    \sum^\infty _{m=0}  \frac{ (ds)^{ m} }{ m !}   |k_i|^m  \|\Box_k f\|_p  \lesssim   \sum^\infty _{m=0}  \frac{  (10ds)^{ m} }{ m !}      \| f\|_p \le M e^{10ds}.
\end{align}
 If $|k_i| > 10$,  in view of Nikol'skij's inequality,  for some $L>(d/(1\wedge p) -1/2)$,
\begin{align} \label{3.14}
    |k_i|^m  \|\Box_k f\|_p & \lesssim  |k_i|^m \sum_{|l|_\infty\le 1} \|\sigma_{k+l} |\xi_i|^{-m} \|_{M_p} \|\Box_k \partial^{m}_{x_i} f\|_p \nonumber\\
      & \lesssim   2^m (1 + m \langle k_i\rangle^{-1}+...+  m(m+1)...(m+L-1) \langle k_i\rangle^{-L})  \|  \partial^{m}_{x_i} f\|_p \nonumber\\
       & \lesssim   2^m  (1+ (m+L)^{L})  \|  \partial^{m}_{x_i} f\|_p .
\end{align}
By \eqref{3.12} and \eqref{3.14},  if $f\in G_{1, p}$ and $|k_i| > 10$,   then for any $s<\rho /2d$,
\begin{align} \label{3.15}
    \sum^\infty _{m=0}  \frac{  (ds)^{ m} }{ m !}  |k_i |^m  \|\Box_k f\|_p  \lesssim \sum^\infty _{m=0} \frac{(2ds)^{m}}{\rho^m } (1+(m+L)^{L})    \lesssim 1.
\end{align}
By \eqref{3.13} and \eqref{3.15},  we have $f\in E^s_{p, \infty}$. It follows that $ G_{1, p}  \subset \bigcup_{s>0} E^s_{p,  \infty}$.  Noticing that $E^s_{p,  \infty} \subset E^{s/2}_{p,  q}$,  we immediately have $ G_{1, p}  \subset \bigcup_{s>0} E^s_{p,  q}$. $\hfill\Box$


\footnotesize
\begin{thebibliography}{100}

\bibitem{BaBiTa12} H. Bae, A. Biswas and E. Tadmor,  Analyticity and decay estimates
of the Navier--Stokes equations in critical
Besov spaces, Arch. Rational Mech. Anal., {\bf 205} (2012), 963--991.


\bibitem{BGOR} A. B\'{e}nyi,
K. Gr\"ochenig, K.A. Okoudjou and L.G. Rogers, Unimodular Fourier
multiplier for modulation spaces, J. Funct. Anal.,  \textbf{246} (2007), 366--384.

\bibitem{BL76} J. Bergh and J. L\"{o}fstr\"{o}m,
 Interpolation Spaces,   Springer--Verlag,  1976.

\bibitem{Bj66} G. Bj\"orck,  Linear partial differential operators and generalized distributions. Ark. Mat. \textbf{6} (1966) 351--407.

\bibitem{BoDeJa03} A. Bonami, B. Demange, P. Jaming, Hermite functions and uncertainty principles for the Fourier and the windowed Fourier transforms,  Rev. Mat. Iberoamericana \textbf{19} (2003), 23--55.


\bibitem{BoPa08} J. Bourgain and N. Pavlovic,  Ill-posedness of the Navier--Stokes equations in a critical
space in 3D. J. Funct. Anal., $\mathbf{255}$ (2008), 2233--2247.

\bibitem{De73} N. G. de Bruijn,
A theory of generalized functions, with applications to Wigner
distribution and Weyl correspondence,  Nieuw Arch. Wisk. (3),  \textbf{21} (1973), 205-280.

\bibitem{Can97} M. Cannone, A generalization of a theorem by Kato on Navier--Stokes equations, Rev. Mat. Iberoamericana, 13 (1997) 515--541.

\bibitem{CaPl96} M. Cannone and F. Planchon,  Self-similar solutions for Navier--Stokes equations in $\mathbb{R}^3$,  Comm. Partial Differ. Eq., $\mathbf{21}$ (1996), 179--193.

\bibitem{CaTo17} M. Cappiello, J. Toft, Pseudo-differential operators in a Gelfand-Shilov setting,   Math. Nachr. \textbf{290} (2017),  738--755.

\bibitem{ChSh11} A. Cheskidov, and R. Shvydkoy, The regularity of weak solutions of the 3D Navier--Stokes equations in $B^{-1}_{\infty, \infty}$,  Arch. Rational Mech. Anal., {\bf 195} (2011), 159--169.

\bibitem{ChWaWaWo18}    M. Chen, B. Wang,   S. Wang, M. W. Wong,  On dissipative nonlinear evolutional pseudo-differential equations, Appl. Comput. Harmonic Anal., https://doi.org/10.1016/j.acha.2018.04.003


\bibitem{CH} J. Y. Chemin,  { \rm Perfect incompressible fluids,  } Oxford University Press, Oxford,
(1998).

\bibitem{Che99} J. Y. Chemin,  Th\'{e}or\`{e}mes d'unicit\'{e} pour le syst\`{e}me de Navier--Stokes tridimensionnel,
J. Anal. Math., {\bf 77} (1999), 27--50.



\bibitem{ChFaSu12} J. Chen, D. Fan and L. Sun, Asymptotic estimates for unimodular Fourier multipliers on modulation spaces, Discrete Contin. Dyn. Syst.,  \textbf{32}  (2012),   467--485.

\bibitem{Ch99} J. Cho,  A characterization of Gelfand--Shilov space based on Wigner distribution, Commun. Korean Math. Soc., \textbf{14} (1999), 761-767.

    \bibitem{ChChKi96} J. Chung, S. Y. Chung, and D. Kim, Characterization of the Gelfand--Shilov spaces via Fourier transforms, Proc. Amer.Math. Soc. \textbf{124} (1996), 2101--2108.




\bibitem{CoFe78} A.C\'{o}rdoba, C. Fefferman, Wave packets and Fourier integral operators. Commun. Partial
Differ. Equations \textbf{3} (1978), 979--1005.

\bibitem{CorNik08} E. Cordero, F. Nicola,  Some new Strichartz estimates for the Schr\"odinger equation. J. Differential Equations, \textbf{245} (2008),   1945--1974.

 \bibitem{CorNik092} E. Cordero, F. Nicola,   Sharpness of some properties of Wiener amalgam and modulation spaces. Bull. Aust. Math. Soc., \textbf{80} (2009),  105--116.




\bibitem{DoLi09} H. Dong and D. Li, Optimal local smoothing and analyticity rate estimates for the generalized
Navier--Stokes equations. Commun. Math. Sci., {\bf 7} (2009), 67--80.

\bibitem{EsSeSv03} L. Escauriaza, G. Serigin and V. Sverak, $L_{3,\infty}$ solutions of Navier--Stokes equations
and backward uniquness, Uspekhi Mat. Nauk., \textbf{58} (2003), 3--44.

\bibitem{Fei83} H. G. Feichtinger, Modulation spaces on locally
compact Abelian group, Technical Report, University of Vienna, 1983.
Published in: ``Proc. Internat. Conf. on Wavelet and Applications",
99--140.  New Delhi Allied Publishers, India, 2003.


\bibitem{FeGr88} H. G. Feichtinger,  K. Gr\"ochenig,
A unified approach to atomic decompositions via integrable group
representations, Lect. Notes in Math.,  \textbf{1302} (1988), 52--73.

\bibitem{FeWe06} H. G. Feichtinger, F. Weisz, Inversion formulas for the short-time Fourier transform, J. Geom. Anal., \textbf{16} (2006), 507--521.

\bibitem{Fo97} C. Foias, What do the Navier--Stokes equations tell us about turbulence? Harmonic
Analysis and Nonlinear Differential Equations (Riverside, 1995). Contemp. Math., \textbf{208} (1997),
151--180.

\bibitem{FoTe89}  C. Foias and R. Temam,  Gevrey class regularity for the solutions of the Navier--Stokes
equations. J. Funct. Anal., \textbf{87} (1989), 359--369.

\bibitem{GeSh68} I. M. Gelfand and G. E. Shilov, Generalized Functions, Spaces of Fundamental and Generalized Functions,  Vol. 2, Academic Press, New York-London, 1968.

\bibitem{Ge08} P. Germain,  The second iterate for the Navier-Stokes equation. J. Funct. Anal., $\mathbf{255}$ (2008), 2248--2264.

\bibitem{GePaSt07} P. Germain,  N. Pavlovic and G. Staffilani,  Regularity of solutions to the Navier-
Stokes equations evolving from small data in $BMO^{-1}$. Int. Math, Res. Not., $\mathbf{21}$ (2007), 35
pages.

\bibitem{Gr92} P. Gr\"obner, {\it Banachr\"aume Glatter Funktionen und
Zerlegungsmethoden}, Doctoral thesis, University of Vienna, 1992.

\bibitem{Groch01} K. Gr\"ochenig,  {\it Foundations of Time--Frequency
Analysis}, Birkh\"auser, Boston, MA, 2001.

\bibitem{GrZi01}
 K. Gr\"ochenig, G. Zimmermann,   Hardy's theorem and the short-time Fourier transform of Schwartz functions,
 J. London Math. Soc. \textbf{63}   (2001),  205--214.

\bibitem{GrZi04}
K. Gr\"ochenig, G. Zimmermann, Spaces of test functions via the STFT,
J. Funct. Spaces Appl. \textbf{2} (2004),  25--53.





\bibitem{GrKu98} Z. Gruji\v{c} and I. Kukavica, Space analyticity for the Navier--Stokes and related equations with initial data in $L^p$, J. Func. Anal., $\mathbf{152}$ (1998),
  247--466.

\bibitem{HaMoOzWa11} H.  Hajaiej, L. Molinet,T. Ozawa, B. Wang, Necessary and sufficient conditions for the fractional Gagliardo-Nirenberg inequalities and applications to Navier-Stokes and generalized boson equations. Harmonic analysis and nonlinear partial differential equations, 159--175, RIMS K\^{o}ky\^{u}roku Bessatsu, B26, Res. Inst. Math. Sci. (RIMS), Kyoto, 2011.

\bibitem{HaWa14}    J. Han, B. Wang, $\alpha$-modulation spaces (I) scaling, embedding and algebraic properties. J. Math. Soc. Japan \textbf{66} (2014), 1315--1373.

\bibitem{Ho60} L. H$\ddot {\rm o}$rmander,  { \rm Estimates for translation invariant
operators in $L^p$ spaces,  }   Acta Math.,  {\bf 104} (1960),
93--139.



\bibitem{Iw10} T. Iwabuchi, Navier--Stokes equations and nonlinear heat equations
in modulation spaces with negative derivative indices, J. Differential Equations, {\bf 248} (2010), 1972--2002.

\bibitem{Ja76} A. J. E. M. Janssen,
Generalized stochastic processes, T. H. Report, No. 76-WSK-07, Department of Mathematics, Technische Hogeschool Eindhoven, Eindhoven, 1976.

\bibitem{KaKoIt14} K. Kato, M. Kobayashi and S. Ito, Estimates on modulation spaces for Schr\"odinger evolution operators
with quadratic and sub-quadratic potentials. J. Funct. Anal., \textbf{266} (2014), 733--753.

\bibitem{Ka84} T. Kato, Strong $L^p$ solutions of the Navier--Stokes equations in $ \mathbb{{R}}^m$ with applications
to weak solutions,  Math. Z., \textbf{187} (1984), 471--480.

\bibitem{KoTa01} H. Koch, D. Tataru, Well-posedness for the Navier--Stokes equations, Adv. Math., {\rm 157} (2001) 22--35.

\bibitem{KoPi16} S. Kostadinova,  S. Pilipovi\'{c}, K. Sanevac, J. Vindasd, The short-time Fourier transform of distributions of exponential
type and tauberian theorems for shift-asymptotics,  Filomat 30:11 (2016), 3047--3061.

\bibitem{LeRi00}   P.G. Lemari\'{e}-Rieusset, {\rm On the analyticity of mild solutions
for the Navier--Stokes equations}, C. R. Acad. Sci. Paris,  Ser I,
{\bf 330} (2000), 183--186.

\bibitem{LeRi02} P.G. Lemari\'{e}-Rieusset,  Recent Developments in the Navier--Stokes Problem, Chapman
\& Hall/CRC Research Notes in Mathematics, vol \textbf{431},  Boca Raton, 2002.

\bibitem{Li} J.-L. Lions,  { \rm Quelques methods de resolution des probl$\grave{e}$mes
aux limites non lin$\acute{e}$aires,  }  (French)Paris:
Dunod/Gauthier-Villars, 1969.


\bibitem{MiSa06} H. Miura and O. Sawada, On the regularizing rate estimates of Koch-Tataru's solution
to the Navier--Stokes equations,  Asymptot. Anal., \textbf{49} (2006), 1--15.

\bibitem{Ob06} H. M. Obiedat,  A topological characterization of the Beurling-Bj\"ork space $\mathscr{C}_w$ using the short-time Fourier transform,
 Cubo \textbf{8} (2006), 33-45.

\bibitem{Pl96} F. Planchon, Global strong solutions in Sobolev or Lebesgue spaces to the incompressible Navier--Stokes equations in $\mathbb{R}^3$,
Ann. Inst. H. Poincare, AN, {\bf 13} (1996), 319--336.


\bibitem{Pi88} S. Pilipovi\'{c}, Tempered ultradistributions, Boll. Un. Mat. Ital. B (7) \textbf{2} (1988), 235--251.

\bibitem{Si67} J. Sebasti\~{a}o Silva, Les s\'{e}ries de multip\^{o}les des physiciens et la th\'{e}orie des ultradistributions, Math. Ann. \textbf{174} (1967), 109--142.


\bibitem{SuTo07} M. Sugimoto amd N. Tomita, The dilation property of
modulation spaces and their inclusion relation with Besov spaces, J.
Funct. Anal., {\bf 248} (2007), 79--106.

\bibitem{Te06} N. Teofanov,
Ultradistributions and Time--Frequency Analysis,  Operator Theory: Advances and
Applications, Pseudo-Differential Operators
and Related Topics,  Vol. \textbf{164} (2006) 173-191.

\bibitem{Te06a} N. Teofanov,
Modulation spaces, Gelfand-Shilov spaces and pseudodifferential operators,
Sampl. Theory Signal Image Process.,  \textbf{5}  (2006), 225-242

\bibitem{To12} J. Toft, The Bargmann transform on modulation and Gelfand-Shilov spaces, with applications to Toeplitz and pseudo-differential operators, J. Pseudo-Differ. Oper. Appl. \textbf{3} (2012), 145--227.

\bibitem{Tr83} H. Triebel,   Theory of Function Spaces,
  Birkh\"{a}user--Verlag,  1983.

  \bibitem{Wang02} B. Wang, {\rm The limit behavoir of solutions for the
complex Ginzburg--Landau  equation,} Commun. Pure Appl. Math.,
{\bf 55} (2002),  481--508.
\bibitem{Wang04} B. Wang, {\rm Exponential Besov spaces and their
applications to certain evolution equations with dissipation,}
Comm. Pure Appl. Anal., {\bf 3} (2004), 883--919.

\bibitem{Wa15}
            \newblock B.  Wang,
             \newblock Ill-posedness for the Navier-Stokes equation in critical Besov spaces $\dot B^{-1}_{\infty, q}$,
             \newblock {Adv. in Math.}, \textbf{268} (2015), 350--372.

\bibitem{WaHu07} B.   Wang and C.  Huang, Frequency-uniform decomposition method for the generalized BO, KdV and NLS equations,  J. Differential Equations, {\bf 239} (2007), 213--250.

\bibitem{WaHuHaGu11} B.  Wang, Z.   Huo, C.  Hao and Z.  Guo, {\it Harmonic Analysis Methods for Nonlinear Evolution Equations,} World Scientific, 2011.

\bibitem{Wa06}  B.  Wang, L.  Zhao, B. Guo,   {\rm  Isometric decomposition operators, function spaces $E^\lambda_{p,q}$ and their applications to nonlinear evolution equations,}  J. Funct. Anal., {\bf 233} (2006), 1--39.



\bibitem{Yo10} T. Yoneda,  Ill-posedness of the 3D Navier--Stokes equations in a generalized Besov
space near $BMO^{-1}$. J. Funct. Anal., $\mathbf{258}$ (2010), 3376--3387.







\end{thebibliography}
\end{document}